\documentclass[11pt]{article}

\usepackage{latexsym,amsmath,amsfonts}

\begin{document}
\renewcommand\arraystretch{1.8}
\def\no{\noindent}
\def\eps{\varepsilon}
\def\na{\nabla}
\def\un{\u\cdot\na}
\def\a{\alpha}
\def\ue{u^\eps}
\def\be{b^\eps}
\def\ppe{p^\eps}
\def\r{ {\mathbb{R}}^2}
\def\t{ {\mathbb{T}}^2}
\def\ber{\begin{eqnarray}}
\def\er{\end{eqnarray}}
\def\pt{\partial_t}
\def\pp{p^{\eps}}
\def\ppp{p^{\eps_1,0}}
\def\p{\phi^{\eps}}
\def\b{b^{\eps}}
\def\bb{b^{\eps_1,0}}
\def\cn{\cdot\nabla}
\def\u{u^{\eps}}
\def\ubl{{\underline{u}}^{\eps}}
\def\ubr{{\bar{u}}^{\eps}}
\def\bbl{{\underline{b}}^{\eps}}
\def\bbr{{\bar{b}}^{\eps}}
\def\ube{u_B^{\eps}}
\def\bbe{b_B^{\eps}}
\def\urr{u_R^{\eps}}
\def\brr{b_R^{\eps}}
\def\ubel{u_{B+}^{\eps}}
\def\uber{u_{B-}^{\eps}}
\def\bbel{b_{B+}^{\eps}}
\def\bber{b_{B-}^{\eps}}
\def\ubet{u_{B3}^{\eps}}
\def\bbet{b_{B3}^{\eps}}
\def\urrt{u_{R3}^{\eps}}
\def\brrt{b_{R3}^{\eps}}
\def\uble{{\underline{u}}^{\eps}}
\def\ubre{{\bar{u}}^{\eps}}
\def\bble{{\underline{b}}^{\eps}}
\def\bbre{{\bar{b}}^{\eps}}
\def\ubee{u_B^{\eps}}
\def\bbee{b_B^{\eps}}
\def\urre{u_R^{\eps}}
\def\brre{b_R^{\eps}}
\def\ubele{u_{B+}^{\eps}}
\def\ubere{u_{B-}^{\eps}}
\def\bbele{b_{B+}^{\eps}}
\def\bbere{b_{B-}^{\eps}}
\def\ubete{u_{B3}^{\eps}}
\def\bbete{b_{B3}^{\eps}}
\def\urrte{u_{R3}^{\eps}}
\def\brrte{b_{R3}^{\eps}}
\def\v{v^{\eps}}
\def\ut{\tilde{u}^{\eps}}
\def\vt{\tilde{v}^{\eps}}
\def\uu{u^{\eps_1,0}}
\def\bu{\underline{b}}
\def\ub{\bu}
\def\up{{\underline{\phi}}^{\eps}}
\def\bp{{\bar\phi}^{\eps}}
\def\br{\bar\rho}
\def\ur{\underline\rho}
\def\w{w^\eps}
\def\lw{{\tilde w}^{\eps}_i}
\def\tw{{\tilde w}^{\eps}}
\def\tbw{{\tilde {w}}^{\eps}}
\def\tu{\tilde u}
\def\wl{\|\w_i\|_{L^2}}
\def\n{\nonumber}
\def\i{\int_\omega}
\def\il{\|_{L^2}}
\def\tp{\tilde{\phi}^{\eps}}
\def\tup{\tilde{\underline{\phi}}^{\eps}}
\def\tbp{\tilde{\bar{\phi}}^{\eps}}
\def\uh{\|\tu^0\|_{L^\infty}}
\def\cnx{\cdot\nabla_{x, y}}

\newtheorem{lemma}{\bf Lemma}[section]

\numberwithin{equation}{section}

       \newtheorem{theorem}{\bf Theorem}[section]

       \newtheorem{prop}[theorem]{\bf Proposition}

       \newtheorem{defi}[theorem]{\bf Definition}

       \newtheorem{remark}[lemma]{\bf Remark}

       \newtheorem{example}{\bf Example}

      %

%


\begin{center}
{\LARGE Boundary Layer Problems in the Viscosity-Diffusion Vanishing
Limits for the Incompressible MHD Systems} \\
\vspace{0.4cm} {\bf Shu Wang$^a$, Zhouping Xin$^b$}\\
{$^a$\it College of Applied Sciences, Beijing University of Technology, Ping Le Yuan 100\\
Beijing 100124, P. R. China\\
Email: wangshu@bjut.edu.cn}\\
{$^b$ \it The Institute of Mathematical Sciences, The Chinese University of Hong Kong\\
satin, N. T., Hong Kong\\
Email: zpxin@ims.cuhk.edu.hk}
\end{center}
\begin{center}
\begin{minipage}{12cm}
{\bf Abstract:}In this paper, we we study boundary layer problems
for the incompressible MHD systems in the presence of physical
boundaries with the standard Dirichlet boundary conditions with
small generic viscosity and diffusion coefficients. We identify a
non-trivial class of initial data for which we can establish the
uniform stability of the Prandtl's type boundary layers and prove
rigorously that the solutions to the viscous and diffusive
incompressible MHD systems converges strongly to the superposition
of the solution to the ideal MHD systems with a Prandtl's type
boundary layer corrector. One of the main difficulties is to deal
with the effect of the difference between viscosity and diffusion
coefficients and to control the singular boundary layers resulting
from the Dirichlet boundary conditions for both the viscosity and
the magnetic fields. One key derivation here is that for the class
of initial data we identify here, there exist cancelations between
the boundary layers of the velocity field and that of the magnetic
fields so that one can use an elaborate energy method to take
advantage this special structure. In addition, in the case of fixed
positive viscosity, we also establish the stability of diffusive
boundary layer for the magnetic field and convergence of solutions
in the limit of zero magnetic diffusion for general initial data.


\end{minipage}
\end{center}
\date{July 15, 2016}

\bigskip

\newpage

\section{Introduction}
We consider in this paper boundary layer problems, zero
viscosity-diffusion vanishing inviscid limit and zero magnetic
diffusion vanishing limit for the three/two-dimensional
incompressible viscous and diffusive magnetohydrodynamic(MHD)
systems with Dirichlet boundary (no-slip characteristic) boundary
conditions \ber &&\pt\ue+\ue\cn\ue+\nabla \ppe-\eps_1 \Delta\ue
=\be\cn\be, \mbox{ in } \Omega\times(0, T),\label{vmhde1}\\
&&\pt\be+\ue\cn\be-\eps_2 \Delta \be=\be\cn\ue,\quad \mbox{ in }
\Omega\times(0,
T),\label{vmhde2}\\
&&div\ue=0, div\be=0, \quad\mbox{ in } \Omega\times(0,
T),\label{vmhde3} \\
&&\ue=0, \be=0, \quad\mbox{ in } \partial\Omega\times(0, T),
\label{vmhde4}\\
&&\ue(t=0)=\ue_0, \be(t=0)=\be_0, \mbox{ with } \quad
div\ue_0=div\be_0=0 \quad\mbox{ on } \quad\Omega , \label{vmhde5}\er
where $\Omega=\omega\times [0, h]$ or $\Omega=\omega\times (0,
\infty)$, and $\omega=\mathbb{T}^2$ or $\mathbb{R}^2$ in three
dimensional case for MHD systems, and $\omega=\mathbb{T}^1$ or
$\mathbb{R}^1$ in two dimensional case for MHD systems.
$\Delta=\frac{\partial^2}{\partial x^2}+\frac{\partial^2}{\partial
y^2}+\frac{\partial^2}{\partial z^2}$ or
$\Delta=\frac{\partial^2}{\partial x^2}+\frac{\partial^2}{\partial
z^2}$ is the three-dimensional or two-dimensional Laplace operator,
$\eps_1=\eps_1(\eps)>0$ is the viscosity coefficient and
$\eps_2=\eps_2(\eps)>0$ is the magnetic diffusion coefficient. The
unknown functions $\ue, \ppe, \be$ are the velocity, the pressure
and the magnetic field of MHD.

The well-posedness, regularity and asymptotic limit problem on the
incompressible viscous and diffusive MHD systems
\eqref{vmhde1}-\eqref{vmhde3} in the whole space or with
slip/no-slip boundary conditions have been studied extensively, see
\cite{Bis, CW, CMZ, DL, HX, HXY, ST, Wu1, Wu2, WW1, XXW} and the
references therein. When $\eps_1>0$ and $\eps_2>0$, MHD systems in
the whole space and in the bounded domains with slip boundary
conditions for the velocity and with no-slip boundary condition for
the magnetic field has a unique global classical solution for smooth
initial data when space dimension $d=2$ and has a global weak
solution for a class of initial data when $d=3$, see \cite{DL, ST}.
Some regularity criterions are also given by some authors, see
\cite{CW, CMZ, HX, Wu1, Wu2} and therein references. Cao and Wu
\cite{CW} obtain global regularity for the 2D MHD equations with
mixed partial dissipation and magnetic diffusion for any smooth
initial data. Wu \cite{Wu1} considers the inviscid limit problem of
incompressible viscous MHD systems in the whole space. Xiao, Xin and
Wu \cite{XXW} investigate the solvability, regularity and vanishing
viscosity limit of incompressible viscous MHD systems with slip
without friction boundary conditions.

It should be noted that, as in the zero-viscosity vanishing limit
for the Navier-Stokes equations, see\cite{Al, CW, E, EE, GM, Ka1,
La, Li, MY, M1, M2, SC1, SC2, TW2, XX1} and related references, the
zero-viscosity and diffusion vanishing limit for incompressible
viscous and diffusive MHD system in a bounded domain, with Dirichlet
boundary conditions at the boundary, is a challenging problem due to
the possible appearance of boundary layers. Recently, \cite{WLW}
considers boundary layer problems and zero viscosity-diffusion
vanishing limit of the incompressible MHD system with no-slip
boundary conditions and the case $\eps_1=\eps_2$ or the case of the
the different horizontal and vertical viscosities and magnetic
diffusions. In this paper, we consider the boundary layer problems
for the more general case $\eps_1\not=\eps_2$ and also consider the
zero magnetic diffusion limit.

Letting $\eps_1\to 0 $ and $\eps_2\to 0$ in
\eqref{vmhde1}-\eqref{vmhde5}, one obtains formally the following
inviscid MHD systems \ber &&\pt   u^{0}+  u^{0}\cn u^{0}+\nabla
p^{0}=  b^{0}\cn   b^{0},\quad \mbox{ in } \quad\Omega\times(0, T),\label{imhd1}\\
&&\pt   b^{0}+  u^{0}\cn   b^{0}=  b^{0}\cn u^{0},\quad
\mbox{ in } \quad\Omega\times(0,
T),\label{imhd2}\\
&&div   u^{0}=0, div   b^{0}=0, \quad\mbox{ in }\quad
\Omega\times(0,
T),\label{imhd3} \\
&&  u^{0}\cdot n=\pm u_3^{0}=0,   b^{0}\cdot n=\pm
b_3^{0}=0, \quad\mbox{ in }\quad
\partial\Omega\times(0, T),
\label{imhd4}\\
&&  u^{0}(t=0)=  u^{0}_0,   b^{0}(t=0)=  b^{0}_0, \quad\mbox{ with }
\quad div   u^{0}_0=div   b^{0}_0=0 \quad\mbox{ on } \quad\Omega
.\label{imhd5}\er On the other hand, setting the magnetic diffusion
coefficient $\eps_2\to 0$ in \eqref{vmhde1}-\eqref{vmhde5}, one gets
formally the following viscous MHD systems \ber
&&\pt\uu+\uu\cn\uu-\eps_1\Delta \uu+\nabla
\ppp=\bb\cn\bb,\quad \mbox{ in } \Omega\times(0, T),\label{hvmhd1}\\
&&\pt\bb+\uu\cn\bb=\bb\cn\uu,\quad \mbox{ in } \quad\Omega\times(0,
T),\label{hvmhd2}\\
&&div\uu=0, div\bb=0, \quad\mbox{ in }\quad \Omega\times(0,
T),\label{hvmhd3} \\
&&\uu=0, \bb_3=0, \quad\mbox{ in }\quad \partial\Omega\times(0, T),
\label{hvmhd4}\\
&&\uu(t=0)=\uu_0, \bb(t=0)=\bb_0, \quad\mbox{ with } \quad
div\uu_0=div\bb_0=0 \quad\mbox{ on } \Omega. \label{hvmhd5}\er The
purpose of this paper is to prove rigorously the above formal limits
under some assumptions on initial data or viscosity and diffusive
coefficients. Since we can recover the incompressible Navier-Stokes
equations by taking $\be=0$ in MHD systems
\eqref{vmhde1}-\eqref{vmhde5}, and, hence, the basic difficulties
caused by such as the well-posedness of the Prandtl's type boundary
layer equations, the thickness of the boundary layer and nonlocal
pressure when dealing with the boundary layer problem and zero
viscosity vanishing limit of the incompressible NS equations in
domains with boundaries are kept here. The key ingredients here are
that we will be able to identify a non-trivial class of initial data
for which there exist cancelations between boundary layers of the
velocity field and that of the magnetic field, which make the
stability of the boundary layers and uniform convergence possible.
Hence the method used here can not be extend directly to deal with
the boundary layer problem for incompressible Navier-Stokes
equations. Moreover, the zero magnetic diffusion limit for MHD
systems in a domain with the boundary is also non-trivial due to the
nonlinear coupling of the velocity field and the magnetic field. The
boundary layer problem for the magnetic field can not also be
obtained and it is remained open in general, see Remark 2.4.

The paper is organized as follows. In section 2 we give the main
results of this paper. Section 3 is devoted to the proofs of Main
results, including the constructs of the approximating boundary
layer functions.

\section{The main results}
In this section, we state our main Theorems. For this, we first
recall  the following classical results  on  the  existence of
sufficiently regular solutions to the incompressible ideal MHD
system (see \cite{DL, ST}).

\begin{prop}\label{prop2.1}
Let $(u^{0}_0,   b^{0}_0)$  satisfy  $u_0^{0}, b_0^{0}\in
H^s(\Omega)$, $s>\frac 32+1$, ${div }{u}_0^{0}={div }{b}_0^{0}=0$
and $u_0^0\cdot n|_{\partial\Omega}=b_0^0\cdot
n|_{\partial\Omega}=0$. Then there exist $0<T_*\le \infty$, the
maximal existence time, and a unique smooth solution $({u}^{0},
p^{0}, b^{0})$, also denoted by $({u}^{0,0}, p^{0,0}, b^{0,0})$
below, of the incompressible ideal MHD equations
\eqref{imhd1}-\eqref{imhd5} on $[0,T_*)$ satisfying, for any
$T<T_*$, \ber \sup_{0\le t\le T}\Big(\|(  u^{0}, b^{0})\|_{H^s}
+\|(\pt{ u}^{0}, \pt
  b^{0})\|_{H^{s-1}(\Omega)}\Big)\le C(T).\nonumber \er
Moreover, if $(u^{0}_0, b^{0}_0)$ satisfies $u^{0}_0=\pm b^{0}_0$,
then there is a unique smooth solution $({u}^{0}, p^{0}, b^{0})$ of
incompressible inviscid MHD system satisfying ${u}^{0}(x,y,z,t)=\pm
b^{0}(x,y,z,t)$ $(={u}^{0}_0(x,y,z)=\pm b^{0}_0(x,y,z))$,
$p^0(x,y,z,t)=0.$ Also, there exist the smooth solutions to the
initial boundary problems for three/two dimensional incompressible
ideal MHD systems for the smooth initial data, which maybe not
belong to Sobolev space $H^s$ in unbounded domain, for example, the
shear flow.
\end{prop}
Similarly, for the incompressible MHD with the viscosity
\eqref{hvmhd1}-\eqref{hvmhd5}, it is easy to get the following
result on the existence of sufficiently regular solutions.
\begin{prop}\label{prop2.2} Assume that $\eps_1>0$ be fixed.
Let $(\uu_0, \bb_0)$  satisfy  $\uu_0, \bb_0\in  H^s(\Omega)$,
$s>\frac 32+2$, ${div }\uu_0={div }\bb_0=0$. Then there exist
$0<T_*\le \infty$, the maximal existence time, and a unique smooth
solution $(\uu, \ppp, \bb)$ of the incompressible MHD equations
\eqref{hvmhd1}-\eqref{hvmhd5} on $[0,T_*)$ satisfying, for any
$T<T_*$, \ber \sup_{0\le t\le T}\Big(\|(\uu, \bb)\|_{H^s(\Omega)}
+\|(\pt\uu, \pt
\bb)\|_{H^{s-2}(\Omega)}\Big)\nonumber\\
+\eps_1\int_0^T\|\nabla\uu(x,y,z, t)\|^2_{H^s(\Omega)}dt\le C(T)
\nonumber \er for some positive constant $C(T)$ independent of
$\eps_2$. Moreover, if $\bb_0|_{\partial\Omega}=0$, then
$\bb(x,y,z,t)|_{\partial\Omega}=0$.
\end{prop}
Now we can state the main results of this paper.

For the MHD system \eqref{vmhde1}-\eqref{vmhde5}, we have the
following stability result of Prandtl's type boundary layer for a
class of special initial data.
\begin{theorem}{\bf (Stability of the Prandtl boundary Layer)}\label{th1} Let $(u^0, p^0, b^0)$ be the solution to the
incompressible ideal MHD system \eqref{imhd1}-\eqref{imhd5}. Assume
that $(\ue_0, \be_0)$ strongly converges in $L^2(\Omega)$ to
$({u}_0^{0}, {b}_0^{0})$, where $u^{0}_0, b_0^0\in H^s(\Omega)$,
$s>\frac 32+1$ satisfies $div u_0^0=div b_0^0=0$ and $u_{0}^{0}\cdot
n|_{\partial\Omega}=b_{0}^{0}\cdot n|_{\partial\Omega}=0$. Assume
that $u_0^{0}(x, y, z)=b_0^{0}(x, y, z)$ or $u_0^{0}(x, y,
z)=-b_0^{0}(x, y, z)$. Furthermore, assume that $\eps,
\eps_1,\eps_2$ satisfy the following convergence: \ber
\frac{\eps_1+\eps_2}{\sqrt\eps}\to 0, \quad
\frac{(\eps_1-\eps_2)^2}{\sqrt\eps\eps(\eps_1+\eps_2)}\to 0,\quad
\frac{(\eps_1-\eps_2)^2}{\eps(\eps_1+\eps_2)}\le C\min\{\eps_1,
\eps_2\}\label{eps-assu}\er for some constant $C>0$, independent of
$\eps, \eps_1, \eps_2$, as $\eps\to 0, \eps_1\to 0, \eps_2\to 0$.
Then there exists a global Leray-Hopf weak solutions $(\ue, \ppe,
\be)$ of (\ref{vmhde1})-(\ref{vmhde5}) such that \ber (\ue-u^0,
\be-b^0)\to (0, 0) \quad \mbox{ in } \quad L^\infty(0, T;
L^2(\Omega))\er for any $T: 0<T<\infty$, as viscosity coefficient
$\eps_1\to 0$ and diffusion coefficient $\eps_2\to 0$.

Moreover, if  \ber\|(\ue_0-u^{0}_0,
\be_0-b^{0}_0)\|_{L^2(\Omega)}^2\le C\eps^\kappa,
\kappa>1,\label{initial-e-ass}\er then there exists $C(T)$,
independent of $\eps$, such that, for $0\le t\le T<\infty$, \ber
\|(\ue-u^{0}, \be-b^{0})\|_{L^2(\Omega)}^2\le
C(T)(\eps^{\kappa-1}+\eps_1^2+\eps_2^2+\frac{\eps_1+\eps_2}{\sqrt\eps}
+\frac{(\eps_1-\eps_2)^2}{\eps\sqrt\eps(\eps_1+\eps_2)})
.\label{res1}\er Furthermore, we have the following stronger
$L^\infty$ convergence results for the viscous and diffusive
incompressible 2D MHD systems: Assume that $\omega=\mathbb{T}^1$ and
$u_1^0(x, z=0)=b_1^0(x, z=0)=const$ (for example, the
two-dimensional shear flow) and $\eps_1=\eps_2$ or that
$\eps,\eps_1,\eps_2$ satisfy suitable relations(stated below). If,
for some suitably large $\kappa>2$,
\ber\|(\ue_0-u^{0}_0-u^\eps_B(t=0),
\be_0-b^{0}_0-b^\eps_B(t=0))\|_{H^s(\Omega)}^2\le C\eps^\kappa,\quad
s>3,\label{initial-e-assc}\er then there exists $C(T)$, independent
of $\eps$, such that, for $0\le t\le T<\infty$, \ber
&&\|(\ue-u^{0}-u^\eps_B,
\be-b^{0}-b^\eps_B)\|_{L^\infty(\Omega\times (0, T)}\le
C(T)\sqrt\eps \mbox{ if } \eps_1=\eps_2=\eps
\label{equ}\\
&&\|(\ue-u^{0}-u^\eps_B,
\be-b^{0}-b^\eps_B)\|_{L^\infty(\Omega\times (0,
T)}\nonumber\\
&\le&C(T)((\frac{\beta_1(\eps)}{\min\{\eps_1,\eps_2\}})^{\frac14}(\beta_2(\eps))^{\frac14}+(\beta_0(\eps))^{\frac
14}(\frac{\beta_4(\eps)}{\min\{\eps_1,\eps_2\}})^{\frac 14})\to 0
\mbox{ when } \eps\to 0.\label{res1c}\er Here $u^\eps_B, b^\eps_B,
\beta_0(\eps),\beta_1(\eps), \beta_2(\eps), \beta_4(\eps)$ will be
given more precisely in the next section.
\end{theorem}
\begin{remark}\label{remarkwell-e} If $\eps_1=\eps_2$ or $\eps_1=\eps, \eps_2=\eps+\eps^{\alpha+1}$ with $\alpha>\frac 12$,
then the assumption \eqref{eps-assu} holds. The boundary layers for
the velocity field and the magnetic field in Theorem \ref{th1}
occurs, which is the standard Prandtl boundary layer and will be
given in the following section 3.1. The proof in establishing the
stability of the Prandtl boundary layer here depends strongly upon
the special structure of the solution $(u^0, p^0, b^0)$ to the
inviscid MHD system, i.e. $u^0=\pm b^0$, which yields to that there
exists the cancelation between the Prandtl boundary layer of the
velocity and the one of the magnetic field. Also, the proof of
Theorem \ref{thm1} in the case $\eps_1\not=\eps_2$ is more complex
than that of the special case $\eps_1=\eps_2>0$ discussed in the
paper \cite{WLW}, and here we need some new techniques. Of course,
if $u_0^{0}|_{\partial\Omega}=\pm b_0^{0}|_{\partial\Omega}=0$ in
Theorem \ref{th1}, called well-prepared initial data, then no
Prandtl's type boundary layer occurs.\end{remark}
\begin{remark}\label{remarkshearfolw} For the shear flow $(u^0, p^0, b^0)(z,t)=(u^0_1(z),
u^0_2(z),0, 0, b^0_1(z), b_2^0(z), 0)$ of the incompressible
inviscid MHD system, if $u^0=\pm b^0$, then the Prandtl boundary
layer for viscosity and diffusive MHD system is stable by using
Theorem \ref{th1}.
\end{remark} For the incompressible MHD system
(\ref{vmhde1})-(\ref{vmhde5}), we also have the following result on
the zero magnetic diffusion limit when one fix the viscosity
coefficient $\eps_1>0$.
\begin{theorem}{\bf (The Zero Magnetic Diffusion Limit)}\label{thm1.2}
Let $\eps_1>0$ be fixed. Let $\eps=\eps_2\to 0$. Let us assume that
$(\u_0, \b_0)$ strongly converges in $L^2(\Omega)$ to $(\uu_0,
\bb_0)$, where $(\uu_0, \bb_0)\in H^s(\Omega)$, $s>\frac 32+2$.
Assume also that $(\uu, \ppp, \bb)$ is the smooth solution to the
system \eqref{hvmhd1}-\eqref{hvmhd5} defined on $[0, T*)$ with
$0<T*\le \infty$, given by Proposition \ref{prop2.2}. Then there
exists global Leray-Hopf weak solutions $(\u, \pp, \b)$ of
(\ref{vmhde1})-(\ref{vmhde5}) such that \ber (\u, \b)-(\uu, \bb))\to
(0, 0) \quad \mbox{ in } \quad L^\infty(0, T;
L^2(\Omega))\label{conv2}\er for any $T: 0<T<T^*$, as $\eps_2\to 0$.

Moreover, if \ber\|(\u_0-\uu_0, \b_0-\bb_0)\|_{L^2(\Omega)}^2\le
C(\sqrt{\eps_2})^{1-\tau}\label{initialass1}\er for any given
$0\le\tau<1$, then there exists $C(T)$, independent of $\eps_2$,
such that \ber \|(\u, \b)-(u^{\eps_1,0},
b^{\eps_1,0})\|_{L^2(\Omega)}^2\le
C(T)(\sqrt{\eps_2})^{1-\tau}.\label{estrate1}\er
\end{theorem}\begin{remark}This is also one boundary layer problem here.
In fact, if $b^{\eps_1,0}_0|_{\partial\Omega}\not=0$, then there
occurs the boundary layer for the magnetic field due to the
difference of the boundary conditions between $\be$ and
$b^{\eps_1,0}$ in the boundary $\partial\Omega$ of the domain.
\end{remark}
\begin{remark}When one replaces the viscosity term $\eps_1\Delta u^\eps$ by
$\eps_1\partial_z^2u^\eps$ in the system
(\ref{vmhde1})-(\ref{vmhde5}), similar zero magnetic diffusion limit
result in Theorem \ref{thm1.2} as $\eps_2\to 0$ holds. Note that the
non-degeneration of the normal direction of the boundary plays a key
role in establishing the stability of the boundary layer. In fact,
the zero magnetic diffusion limit of the following MHD system \ber
&&\pt\ue+\ue\cn\ue+\nabla \ppe
=\be\cn\be, \mbox{ in } \Omega\times(0, T),\nonumber\\
&&\pt\be+\ue\cn\be-\eps_2 \Delta \be=\be\cn\ue,\quad \mbox{ in }
\Omega\times(0,
T),\nonumber\\
&&div\ue=0, div\be=0, \quad\mbox{ in } \Omega\times(0,
T),\nonumber \\
&&\ue_3=0, \be=0, \quad\mbox{ in } \partial\Omega\times(0, T),
\nonumber\\
&&\ue(t=0)=\ue_0, \be(t=0)=\be_0, \mbox{ with } \quad
div\ue_0=div\be_0=0 \quad\mbox{ on } \quad\Omega \nonumber\er or
\ber &&\pt\ue+\ue\cn\ue+\nabla \ppe-\eps_1 \Delta_{x,y}\ue
=\be\cn\be, \mbox{ in } \Omega\times(0, T),\nonumber\\
&&\pt\be+\ue\cn\be-\eps_2 \Delta \be=\be\cn\ue,\quad \mbox{ in }
\Omega\times(0,
T),\nonumber\\
&&div\ue=0, div\be=0, \quad\mbox{ in } \Omega\times(0,
T),\nonumber \\
&&\ue_3=0, \be=0, \quad\mbox{ in } \partial\Omega\times(0, T),
\nonumber\\
&&\ue(t=0)=\ue_0, \be(t=0)=\be_0, \mbox{ with } \quad
div\ue_0=div\be_0=0 \quad\mbox{ on } \quad\Omega \nonumber\er is
open if $\be|_{\partial\Omega}\not=0$, which yields the appearance
of the boundary layer. Here $\Delta_{x,y}=\frac{\partial^2}{\partial
x^2}+\frac{\partial^2}{\partial y^2}$ and $\eps_2=\eps$.
\end{remark} \begin{remark}One of the main difficulties to establish the zero viscosity and diffusion limit
is to deal with the terms related to the boundary layers in the
error equations. However, the proof is elementary if there occurs no
boundary layers. For example, it is easy to prove that there exists
a $T>0$, independent of $\eps$, such that the solution $(u^\eps,
p^\eps, b^\eps)$ of the system \ber &&\pt\ue+\ue\cn\ue+\nabla
\ppe-\eps_1 \Delta\ue
=\be\cn\be, \mbox{ in } \Omega\times(0, T),\label{nvmhde1}\\
&&\pt\be+\ue\cn\be-\eps_2 \Delta \be=\be\cn\ue,\quad \mbox{ in }
\Omega\times(0,
T),\label{nvmhde2}\\
&&div\ue=0, div\be=0, \quad\mbox{ in } \Omega\times(0,
T),\label{nvmhde3} \\
&&\ue=u^0, \be=b^0, \quad\mbox{ in } \partial\Omega\times(0, T),
\label{nvmhde4}\\
&&\ue(t=0)=\ue_0, \be(t=0)=\be_0, \mbox{ with } \quad
div\ue_0=div\be_0=0 \quad\mbox{ on } \quad\Omega  \label{nvmhde5}\er
converges to the solution $(u^0, p^0, b^0)$ of the ideal MHD system
(\ref{imhd1})-(\ref{imhd5}) in the interval $[0, T]$ in some kinds
of norm, for example, in $L^2(\Omega)$, when $\eps_1\to 0$ and
$\eps_2\to 0$. Here the function $(u^0, b^0)$ given in the boundary
condition \eqref{nvmhde4} of the system
\eqref{nvmhde1}-\eqref{nvmhde5} is determined by the solution of the
ideal MHD system (\ref{imhd1})-(\ref{imhd5}).
\end{remark}
\section{The proofs of Theorems \ref{th1} and \ref{thm1.2}}
We will prove Theorems \ref{th1} and \ref{thm1.2} when
$\Omega=\omega\times [0, h]$ for $0<h<\infty$ and
$\omega=\mathbb{T}^2$. The other cases (for example, when
$h=\infty$) are similar. Our proof is based on the asymptotic
analysis with multiple scales and the classical energy method. We
will divide the proof into two cases. For well-prepared initial data
$b^0|_{\partial\Omega}=b^{0}_0|_{\partial\Omega}=0$ in Theorem
\ref{th1} or $\bb_0|_{z=0}=0$, there is no boundary layer for the
magnetic field, and, hence, there is no boundary layer for the
velocity in the proof of Theorem \ref{th1} for the case of
well-prepared initial data due to $u_0^0(x,y,z)=b_0^0(x,y,z)$. For
the general initial data, i.e.,
$b^{0}_0|_{\partial\Omega}=u^{0}_0|_{\partial\Omega}\not=0$ or
$\bb_0|_{z=0}\not=0$, if one uses energy method to estimate the
error function $(\ue-u^{0}, \be-b^{0})$ or $(\ue-\uu,\be-\bb)$, then
integrations by parts introduce some terms which  are difficult to
control, because $\ue-u^0, \be-b^0$ or $\be-\bb$ do not vanish at
the boundary. So, for general initial data, one needs to construct
the boundary layer correctors which allow one to recover zero
Dirichlet boundary condition. When $0<h<\infty$, we will construct
the left and right boundary layers respectively. When $h=\infty$, we
will construct only the left boundary layer by taking the right
boundary layer to be zero. Note that there is no boundary layer for
the velocity field in the case of Theorem \ref{thm1.2}.
\subsection{The construction of the boundary layers}
We only construct the boundary layer $(u_B^{\eps},b_B^{\eps})$ for
the viscous and diffusive MHD system with $\eps_1\to 0$ and
$\eps_2\to 0$. In the case that $\eps_1>0$ is a fixed given
constant, the structure and its properties of the boundary layer,
denoted also by $b_B^{\eps}$ with $\eps=\eps_2$, are the same as
$b_B^\eps$ by replacing the function $b^0$ by $b^{\eps_1,0}$ and
$\nu_2^*=(\theta\eps_2)^{1+\tau}$ for any $0\le\tau<1$ and for
$\theta>0$ sufficiently small.

Because the velocity $\ue$ and the magnetic field $\be$ satisfy
respectively the zero Dirichlet boundary condition at the boundary
$z=0, h$, but $u^{0}|_{z=0, h}\not=0$ and $b^{0}|_{z=0, h}\not=0$,
we need respectively construct the correctors for $u^0$ and $b^0$ so
as to recover the zero Dirichlet boundary condition for the error
functions.

Now, we introduce the following exact boundary layers for the
velocity field $u^{0}$: \ber
u_B^{\eps}=u_{B+}^{\eps}+u_{B-}^{\eps},\label{ub}\er
\ber && u_{B+}^{\eps}=\left(\begin{matrix}U^{\eps}_{B+}\\
\ubl_{3}\end{matrix}\right)
=\left(\begin{matrix}\ubl_{1}\\ \ubl_{2}\\
\ubl_{3}\end{matrix}\right)\nonumber\\
&=&\left(\begin{matrix} -u^{0}_1(x, y, 0, t)e^{-\frac
z{\sqrt{\nu_1^*}}}(\rho_1(z)-\rho_1'(z)\sqrt{\nu_1^*})
-u^{0}_1(x, y, 0, t)\rho_1'(z)\sqrt{\nu_1^*}\\
-u^{0}_2(x, y, 0, t)e^{-\frac
z{\sqrt{\nu_1^*}}}(\rho_1(z)-\rho_1'(z)\sqrt{\nu_1^*})
-u^{0}_2(x, y, 0, t)\rho_1'(z)\sqrt{\nu_1^*}\\
\sqrt{\nu_1^*}\partial_zu^{0}_3(x, y, 0, t)\rho_1(z)(e^{-\frac
z{\sqrt{\nu_1^*}}}-1)\end{matrix}\right)\label{ub1}\er and
\ber &&u_{B-}^{\eps} =\left(\begin{matrix}U^{\eps}_{B-}\\
\ubr_{3}\end{matrix}\right)
=\left(\begin{matrix}\ubr_{1}\\ \ubr_{2}\\
\ubr_{3}\end{matrix}\right)\nonumber\\
 &=&\left(\begin{matrix}
-u^{0}_1(x, y, h, t) e^{-\frac {h-z}{\sqrt{
\nu_1^*}}}(\rho_2(z)+\rho_2'(z)\sqrt{\nu_1^*})
+u^{0}_1(x, y, h, t)\rho_2'(z)\sqrt{\nu_1^*}\\
-u^{0}_2(x, y, h, t)e^{-\frac {h-z}{\sqrt{
\nu_1^*}}}(\rho_2(z)+\rho_2'(z)\sqrt{\nu_1^*})
+u^{0}_2(x, y, h, t)\rho_2'(z)\sqrt{\nu_1^*}\\
-\sqrt{\nu_1^*}\partial_zu^{0}_3(x, y, h, t)\rho_2(z)(e^{-\frac
{h-z}{\sqrt{\nu_1^*}}}-1)\end{matrix}\right)\label{ub2}\er where
$\rho_1(z)$ and $\rho_2(z)$ satisfy \ber \rho_1(0)=1,
\rho'_1(0)=\rho''_1(0)=0; \rho_1(z)\equiv 0 \mbox{ for } z\ge\frac
h4; 0\le\rho_1(z)\le 1 \mbox{ for } z\in [0, h]\label{r1}\er and
\ber \rho_2 (h)=1, \rho'(h)=\rho''_2(h)=0; \rho_2(z)\equiv 0 \mbox{
for } z\le\frac {3}4h; 0\le\rho_2(z)\le 1 \mbox{ for } z\in [0,
h].\label{r2}\er Here $u_{B+}^{\eps}$ and $u_{B-}^{\eps}$ are the
left boundary layer at $z=0$ and the right boundary layer at $z=h$
for $u^{0}$ respectively. When $h=\infty$, we can take
$u_{B-}^{\eps}=0$.

Similarly, the exact boundary layers for the magnetic field $ b^{0}$
can be exactly as: \ber
b_B^{\eps}=b_{B+}^{\eps}+b_{B-}^{\eps},\label{bb}\er
\ber && b_{B+}^{\eps}=\left(\begin{matrix}B^{\eps}_{B+}\\
\bbl_{3}\end{matrix}\right)
=\left(\begin{matrix}\bbl_{1}\\ \bbl_{2}\\
\bbl_{3}\end{matrix}\right)\nonumber\\
&=&\left(\begin{matrix} -  b^{0}_1(x, y, 0, t)e^{-\frac
z{\sqrt{\nu_2^*}}}(\rho_1(z)-\rho_1'(z)\sqrt{\nu_2^*})
-  b^{0}_1(x, y, 0, t)\rho_1'(z)\sqrt{\nu_2^*}\\
-  b^{0}_2(x, y, 0, t)e^{-\frac z{\sqrt{
\nu_2^*}}}(\rho_1(z)-\rho_1'(z)\sqrt{\nu_2^*})
-  b^{0}_2(x, y, 0, t)\rho_1'(z)\sqrt{\nu_2^*}\\
\sqrt{\nu_2^*}\partial_z  b^{0}_3(x, y, 0, t)\rho_1(z)(e^{-\frac
z{\sqrt{\nu_2^*}}}-1)\end{matrix}\right) \label{bb1}\er and
\ber &&b_{B-}^{\eps} =\left(\begin{matrix}B^{\eps}_{B-}\\
\bbr_{3}\end{matrix}\right)
=\left(\begin{matrix}\bbr_{1}\\ \bbr_{2}\\
\bbr_{3}\end{matrix}\right)\nonumber\\
 &=&\left(\begin{matrix}
-  b^{0}_1(x, y, h, t) e^{-\frac {h-z}{\sqrt{
\nu_2^*}}}(\rho_2(z)+\rho_2'(z)\sqrt{\nu_2^*})
+  b^{0}_1(x, y, h, t)\rho_2'(z)\sqrt{\nu_2^*}\\
-  b^{0}_2(x, y, h, t)e^{-\frac {h-z}{\sqrt{
\nu_2^*}}}(\rho_2(z)+\rho_2'(z)\sqrt{\nu_2^*})
+  b^{0}_2(x, y, h, t)\rho_2'(z)\sqrt{\nu_2^*}\\
-\sqrt{\nu_2^*}\partial_z  b^{0}_3(x, y, h, t)\rho_2(z)(e^{-\frac
{h-z}{\sqrt{\nu_2^*}}}-1)\end{matrix}\right).\label{bb2}\er Here
$b_{B+}^{\eps}$ and $b_{B-}^{\eps}$ are the left boundary layer at
$z=0$ and the right boundary layer at $z=h$ for $b^{0}$
respectively. $\rho_1(z)$ and $\rho_2(z)$ are given by \eqref{r1}
and \eqref{r2}. When $h=\infty$, take $b_{B-}^{\eps}=0$. For
well-prepared initial data, take $b_{B+}^{\eps}=b_{B-}^{\eps}=0$.

It is easy to verify the following properties.
\begin{lemma}\label{lem1} There is a positive constant
$C$, depending upon $\|(u^{0},   b^{0})\|_{H^s(\Omega)}, s>\frac
32+1$, but independent of $\eps, \eps_1$ and $\eps_2$, such that
\ber &&\|(\ube, \bbe)\|_{L^\infty(\Omega)}\le C;
\quad \|(\partial_zu_{B3}^{\eps}, \partial_zb_{B3}^{\eps})\|_{L^\infty(\Omega)}\le C;\nonumber\\
&&\|(U_B^{\eps}, \nabla_{x, y}U_B^{\eps})\|_{L^2(\Omega)}\le
C\sqrt{\sqrt{\nu_1^*}}; \quad \|(B_B^{\eps}, \nabla_{x,
y}B_B^{\eps})\|_{L^2(\Omega)}\le
C\sqrt{\sqrt{\nu_2^*}};\nonumber\\
&&\|\ubet\|_{L^2(\Omega)}\le C{\sqrt{\nu_1^*}}; \quad
\|\bbet\|_{L^2(\Omega)}\le
C{\sqrt{\nu_2^*}};\nonumber\\
&&\|\partial_zU_{B}^{\eps}\|_{L^2(\Omega)}\le \frac
C{\sqrt{\sqrt{\nu_1^*}}}; \quad
\|\partial_zB_{B}^{\eps}\|_{L^2(\Omega)}\le \frac
C{\sqrt{\sqrt{\nu_2^*}}}\nonumber\\
&&\|\partial_z\ubet\|_{L^2(\Omega)}\le C\sqrt{\sqrt{\nu_1^*}}; \quad
\|\partial_z\bbet\|_{L^2(\Omega)}\le
C\sqrt{\sqrt{\nu_2^*}};\nonumber\\
&&\|(z\partial_zU_{B+}^{\eps},
(h-z)\partial_zU_{B-}^{\eps})\|_{L^2(\Omega)}\le
C{\sqrt{\sqrt{\nu_1^*}}}; \quad
\|(z\partial_zB_{B+}^{\eps}, (h-z)\partial_zB_{B-}^{\eps})\|_{L^2(\Omega)}\le C{\sqrt{\sqrt{\nu_2^*}}};\nonumber\\
&&\|(z^2\partial_zU_{B+}^{\eps},
(h-z)^2\partial_zU_{B-}^{\eps}\|_{L^\infty(\Omega)}\le
C\sqrt{\nu_1^*}; \quad \|(z^2\partial_zB_{B+}^{\eps},
z^2\partial_zB_{B-}^{\eps})\|_{L^\infty(\Omega)}\le
C\sqrt{\nu_2^*};\nonumber \er
\end{lemma}
$\nu_1^*, \nu_2^*$ will be chosen later. Here and in what follows we
use $U, B, U_B, B_B, U_R,\cdots$ or $u_{1,2},b_{1,2}$,
$\{u_B\}_{1,2}$, $\{b_B\}_{1,2}$,$\{u_R\}_{1,2},\cdots$ to denote
the first and the second components of $u, b, u_B, b_B, u_R,\cdots$
respectively. Also, $u_3, b_3, u_{B3}$, $b_{B3},\cdots,$ denote the
third components of the vectors $u, b,u_B, b_B, \cdots$.
\subsection{The  proof of Theorem \ref{th1}}
The global existence of Leray-Hopf weak solutions to the three
dimensional dissipative incompressible MHD system and the global
existence of smooth solutions to the two dimensional dissipative MHD
system can be proven as in \cite{DL, ST, XXW} by Galerkin method,
and also as for Navier-Stokes equations. We omit details here.

Let $(\ue, \be)$ be Leray-Hopf weak solution of MHD system
\eqref{vmhde1}-\eqref{vmhde5}. Decompose $(\ue, \be)$ as
$(u^{0}+\ubee+\urre, b^{0}+\bbee+\brre)$. Taking
$\nu_1^*=\nu_2^*=\eps$ in the subsection 3.1 and using the system
\eqref{imhd1}-\eqref{imhd5}, we have \ber
&&\pt\ubee+\pt\urre+\ue\cdot\nabla\urre+u^{0}\cdot\nabla\ubee+\ubee\cdot\nabla\ubee+\urre\cdot\nabla\ubee+\ubee\cdot\nabla
u^{0} \nonumber\\&&+\urre\cdot\nabla u^{0}
-\eps_1\partial_z^2\urre-\eps_1\partial_z^2u^{0}-\eps_1\partial_z^2\ubee-\eps_1\Delta_{x,
y}\urre-\eps_1\Delta_{x, y}u^{0}-\eps_1\Delta_{x, y}\ubee\nonumber\\
&=&-\nabla
(\ppe-p^{0})+\be\cdot\nabla\brre+  b^{0}\cdot\nabla\bbee+\bbee\cdot\nabla\bbee+\brre\cdot\nabla\bbee\nonumber\\
&&+\bbee\cdot\nabla  b^{0}+\brre\cdot\nabla  b^{0}, \quad \mbox{ in }\quad \Omega\times(0, T),\label{mhdre1}\\
&&\pt\bbee+\pt\brre+\ue\cdot\nabla\brre+u^{0}\cdot\nabla\bbee+\ubee\cdot\nabla\bbee+\urre\cdot\nabla\bbee+\ubee\cdot\nabla
b^{0} \nonumber\\&&+\urre\cdot\nabla  b^{0}
-\eps_2\partial_z^2\brre-\eps_2\partial_z^2
b^{0}-\eps_2\partial_z^2\bbee-\eps_2\Delta_{x,
y}\brre-\eps_2\Delta_{x, y}  b^{0}-\eps_2\Delta_{x, y}\bbee\nonumber\\
&=&\be\cdot\nabla\urre+  b^{0}\cdot\nabla\ubee+\bbee\cdot\nabla\ubee+\brre\cdot\nabla\ubee+\bbee\cdot\nabla u^{0}\nonumber\\
&&+\brre\cdot\nabla u^{0}, \quad \mbox{ in }\quad \Omega\times(0, T),\label{mhdre2}\\
&&div\ue=div\be=div\urre=div\brre=0,\mbox{ in } \Omega\times(0, T),\label{mhdre3}\\
&&\urre=\brre=0,  \quad\mbox { on } \quad (x, y)\in \omega, z=0, h,
0\le t\le
T\label{mhdre4}\\
&&\urre(t=0)=\ue(0)-u^{0}(0)-\ubee(t=0),\nonumber\\
&&\brre(t=0)=\be(0)-  b^{0}(0)-\bbee(t=0), \quad \mbox{ on } \quad
\Omega .\label{mhdre5}\er Thanks to the fact that $u^0_0=\pm b^0_0$,
where $(u^0,p^0,b^0)(x,y,z,t)=(u^0_0,0,b^0_0)$ is the special
solution to the incompressible MHD equation, we have that
$u_B^\eps=\pm b_B^\eps$, which shows that there exist cancelations
between the boundary layers of the velocity and the magnetic field,
and, hence, it follows from \eqref{mhdre1} and \eqref{mhdre2} that
\ber
&&\pt(\urre-\brre)+\ue\cdot\nabla(\urre-\brre)-\frac{\eps_1+\eps_2}2\Delta
(\urre-\brre)-\frac{\eps_1-\eps_2}2\Delta (\urre+\brre)
\nonumber\\
&=&-\nabla
(\ppe-p^{0})-\be\cdot\nabla(\urre-\brre)+\frac{\eps_1-\eps_2}2\Delta(u^\eps_B+b_B^\eps)+\frac{\eps_1-\eps_2}2\Delta(u^0+b^0)\label{mhdre-s1}\er
or \ber
&&\pt(\urre+\brre)+\ue\cdot\nabla(\urre+\brre)-\frac{\eps_1+\eps_2}2\Delta
(\urre+\brre)-\frac{\eps_1-\eps_2}2\Delta (\urre-\brre)
\nonumber\\
&=&-\nabla
(\ppe-p^{0})+\be\cdot\nabla(\urre+\brre)+\frac{\eps_1-\eps_2}2\Delta(u^\eps_B-b_B^\eps)+\frac{\eps_1-\eps_2}2\Delta(u^0-b^0)\label{mhdre-s2}\er
Here we has used the relation $\eps_1\Delta u^\eps_R-\eps_2\Delta
b^\eps_R=\frac{\eps_1+\eps_2}2\Delta
(u^\eps_R-b^\eps_R)+\frac{\eps_1-\eps_2}2\Delta(u^\eps_R+b^\eps_R)$
or $\eps_1\Delta u^\eps_R+\eps_2\Delta
b^\eps_R=\frac{\eps_1+\eps_2}2\Delta
(u^\eps_R+b^\eps_R)+\frac{\eps_1-\eps_2}2\Delta(u^\eps_R-b^\eps_R)$.
Noting that there appears the term $-\frac{\eps_1-\eps_2}2\Delta
(\urre+\brre)$ in the system \eqref{mhdre-s1} due to the fact that
$\eps_1\not=\eps_2$, and, hence, one can not adopt the techniques in
\cite{WLW} here, so a new idea is needed to deal with the current
case. Now, using \eqref{mhdre3}-\eqref{mhdre4} and taking the scalar
product of \eqref{mhdre-s1} (or \eqref{mhdre-s2}) with $\urre-\brre$
(or $\urre+\brre$), we get, for the case of $u^0=b^0$, that \ber
&&\frac 12\frac d{dt}\int
|\urre-\brre|^2+\frac{\eps_1+\eps_2}2\int|\nabla(\urre-\brre)|^2\nonumber\\
&=&-\frac{\eps_1-\eps_2}2\int\nabla((u_R^\eps+b_R^\eps)+(u_B^\eps+b_B^\eps)+(u^0+b^0))\cdot\nabla(u_R^\eps-b_R^\eps)\label{es-s1}\er
and for the case of $u^0=-b^0$, that \ber &&\frac 12\frac d{dt}\int
|\urre+\brre|^2+\frac{\eps_1+\eps_2}2\int|\nabla(\urre+\brre)|^2\nonumber\\
&=&-\frac{\eps_1-\eps_2}2\int\nabla((u_R^\eps-b_R^\eps)+(u_B^\eps-b_B^\eps)+(u^0-b^0))\cdot\nabla(u_R^\eps+b_R^\eps),\label{es-s2}
\er which yields to, by using \eqref{mhdre5} and the properties of
the boundary layer functions and with the help of the Young's
inequality and the assumption \eqref{initial-e-ass} on initial data,
that, for $0\le t\le T<\infty$, \ber &&\int
|(\urre-\brre)(t)|^2+(1-\delta)(\eps_1+\eps_2)\int_0^t\int|\nabla(u_R^\eps-b_R^\eps)|^2\nonumber\\
&\le&
\int
|(\urre-\brre)(t=0)|^2+\frac{(\eps_1-\eps_2)^2}{4\delta(\eps_1+\eps_2)}[\int_0^t\int(|\nabla
u_R^\eps|^2+|\nabla b_R^\eps|^2)+Ct+\frac C{\sqrt\eps}t]\nonumber\\
&\le&
C\eps^\kappa+\frac{(\eps_1-\eps_2)^2}{4\delta(\eps_1+\eps_2)}\int_0^t\int(|\nabla
u_R^\eps|^2+|\nabla
b_R^\eps|^2)+C\frac{(\eps_1-\eps_2)^2}{4\delta\sqrt\eps(\eps_1+\eps_2)}\label{est-se1}\er
or \ber &&\int
|(\urre+\brre)(t)|^2+(1-\delta)(\eps_1+\eps_2)\int_0^t\int|\nabla(u_R^\eps+b_R^\eps)|^2\nonumber\\
&\le&
C\eps^\kappa+\frac{(\eps_1-\eps_2)^2}{4\delta(\eps_1+\eps_2)}\int_0^t\int(|\nabla
u_R^\eps|^2+|\nabla
b_R^\eps|^2)+C\frac{(\eps_1-\eps_2)^2}{4\delta\sqrt\eps(\eps_1+\eps_2)}\label{est-se2}\er
for some constant $C>0$ and $\delta>0$ independent of
$\eps,\eps_1,\eps_2$, and $\kappa>1$ .

Now, by taking the scalar product of \eqref{mhdre-s1} with $\urre$
and the scalar product of \eqref{mhdre-s2} with $\brre$, one can get
that \ber \frac 12\frac d{dt}\int
(|\urre|^2+|\brre|^2)+\eps_1\int|\nabla\urre|^2+\eps_2\int|\nabla\brre|^2
=\sum_{i=1}^{13}J_i,\label{est-se3} \er where $J_i, i=1, \cdots,
13,$ are given respectively as follows \ber
&&J_1=-\int\nabla(\ppe-p^0)\urre;\nonumber\\
&&J_2=-\int \pt\ubee \urre-\int\pt\bbee\brre;\nonumber\\
&&J_3=-\int\ue\cdot\nabla\urre\urre-\int\ue\cdot\nabla\brre\brre;\nonumber\er
\ber &&J_4=-\int u^{0}\cdot\nabla\ubee\urre-\int
u^{0}\cdot\nabla\bbee\brre+\int  b^{0}\cdot\nabla\bbee\urre+\int  b^{0}\cdot\nabla\ubee\brre; \nonumber\\
&&J_5=-\int\ubee\cdot\nabla\ubee\urre-\int\ubee\cdot\nabla\bbee\brre+\int\bbee\cdot\nabla\bbee\urre+\int\bbee\cdot\nabla\ubee\brre;\nonumber\\
&&J_6=-\int\urre\cdot\nabla\ubee\urre-\int\urre\cdot\nabla\bbee\brre +\int\brre\cdot\nabla\bbee\urre+\int\brre\cdot\nabla\ubee\brre\nonumber\\
&&J_7=-\int\ubee\cdot\nabla u^{0}\urre-\int\ubee\cdot\nabla  b^{0}\brre+\int\bbee\cdot\nabla  b^{0}\urre+\int\bbee\cdot\nabla u^{0}\brr;\nonumber\\
&&J_8=-\int\urre\cdot\nabla u^{0}\urre-\int\urre\cdot\nabla  b^{0}\brre +\int\brre\cdot\nabla  b^{0}\urre+\int\brre\cdot\nabla u^{0}\brre; \nonumber\\
&&J_{9}=\int\be\cdot\nabla\brre\urre+\int\be\cdot\nabla\urre\brre;\nonumber\er
\ber &&J_{10}=\eps_1\int\partial_z^2\ubee\urre+\eps_2\int\partial_z^2 \bbee\brre; \quad J_{11}=\eps_1\int\partial_z^2u^{0}\urre+\eps_2\int\partial_z^2  b^{0}\brre;\nonumber\\
&&J_{12}=\eps_1\int\Delta_{x, y}\ubee\urre+\eps_2\int\Delta_{x,
y}\bbee\brre; \quad J_{13}=\eps_1\int\Delta_{x,
y}u^{0}\urre+\eps_2\int\Delta_{x, y}b^{0}\brre.\nonumber\er We now
bound each of $J_i, i=1, \cdots, 13$. In the sequel, $C$ denotes any
constant depending only upon $h$. Also, in the following, we just
consider the case of $u^0=b^0$ in Theorem \ref{th1} because the case
of $u^0=-b^0$ can be treated similarly by using \eqref{est-se2}

1) First, using the fact that $u^0=b^0$, which is independent of the
time $t$, and, hence, $u^\eps_B=b^\eps_B$,  we get
$$J_2=J_4=J_5=J_7=0.$$

2) Second, direct computation gives \ber
J_1=-\int\nabla(\ppe-p^0)\urre=\int (\ppe-p^0) div\urre
=0,\label{est-e1}\er \ber
J_3&=&-\int\ue\cdot\nabla\urre\urre-\int\ue\cdot\nabla\brre\brre\nonumber\\
&=&\frac 12\int\ue\cdot\nabla(|\urre|^2)+\frac
12\int\ue\cdot\nabla(|\brre|^2)=0,\label{est-e2}\er To estimate the
term $J_6$, noting that $u^\eps_B=b^\eps_B$ one gets by using the
estimate \eqref{est-se1} that \ber J_6&=&-\int
(\urre-\brre)\cdot\nabla u_B^\eps (\urre+\brre)\nonumber\\
&=&-\int (\urre-\brre)_{1,2}\cdot\nabla_{x,y} (u_{B}^\eps)_{1,2}
(\urre+\brre)_{1,2}-\int
(\urre-\brre)_{1,2}\cdot\nabla_{x,y} (u_{B}^\eps)_{3} (\urre+\brre)_{3}\nonumber\\
&&-\int (\urre-\brre)_{3}\cdot\nabla_{z} (u_{B}^\eps)_{1,2}
(\urre+\brre)_{1,2}-\int (\urre-\brre)_{3}\cdot\nabla_{z}
(u_{B}^\eps)_{3} (\urre+\brre)_{3}\nonumber\\
&\le&
C\int(|\urre|^2+|\brre|^2)+C\int\frac{|(\urre-\brre)_3|^2}\eps\nonumber\\
&\le&C\int(|\urre|^2+|\brre|^2)+C\eps^{\kappa-1}+\frac{(\eps_1-\eps_2)^2}{4\delta\eps(\eps_1+\eps_2)}\int_0^t\int(|\nabla
u_R^\eps|^2+|\nabla
b_R^\eps|^2)\nonumber\\
&&+C\frac{(\eps_1-\eps_2)^2}{4\delta\eps\sqrt\eps(\eps_1+\eps_2)},
\label{est-e3}\er \ber J_8\le 2\max\{\|\nabla
u^0\|_{L^\infty},\|\nabla b^0\|_{L^\infty}\}\int
(|\urre|^2+|\brre|^2),\label{est-e4}\er \ber J_9=\int b^\eps
\cdot\nabla(\urre\cdot\brre)=\int div b^\eps
(\urre\cdot\brre)=0,\label{est-e5}\er \ber
J_{10}&=&-\eps_1\int\partial_z\ubee\partial_z\urre-\eps_2\int\partial_z
\bbee\partial_z\brre\nonumber\\
&\le& \delta\eps_1\int|\nabla u_R^\eps|^2+\delta\eps_2\int|\nabla
b_R^\eps|^2+\frac{C(\eps_1+\eps_2)}{\sqrt\eps},\label{est-e10} \er
\ber J_{11}+J_{12}+J_{13}\le
C\int(|\urre|^2+|\brre|^2)+C(\eps_1^2+\eps_2^2).\label{est-e6}\er
Combining \eqref{est-se3} with \eqref{est-e1}-\eqref{est-e6}
together and using the assumption \eqref{eps-assu} in Theorem
\ref{th1}, we have \ber &&\frac 12 \frac
d{dt}\int(|\urre|^2+|\brre|^2)+(1-\delta)\eps_1\int
|\nabla\urre|^2+(1-\delta)\eps_2\int |\nabla\brre|^2\nonumber\\
&\le& C\int
(|\urre|^2+|\brre|^2)+\frac{(\eps_1-\eps_2)^2}{4\delta\eps(\eps_1+\eps_2)}\int_0^t\int(|\nabla
u_R^\eps|^2+|\nabla
b_R^\eps|^2)\nonumber\\
&&+C(\eps^{\kappa-1}+\eps_1^2+\eps_2^2+\frac{\eps_1+\eps_2}{\sqrt\eps})
+C\frac{(\eps_1-\eps_2)^2}{4\delta\eps\sqrt\eps(\eps_1+\eps_2)}\nonumber\\
&\le& C[\int (|\urre|^2+|\brre|^2)+\eps_1\int_0^t\int|\nabla
u_R^\eps|^2+\eps_2\int_0^t\int|\nabla
b_R^\eps|^2]\nonumber\\
&&+C(\eps^{\kappa-1}+\eps_1^2+\eps_2^2+\frac{\eps_1+\eps_2}{\sqrt\eps})
+C\frac{(\eps_1-\eps_2)^2}{4\delta\eps\sqrt\eps(\eps_1+\eps_2)}\label{est-se4}\er
It follows from \eqref{est-se4} and by Gronwall's inequality that
\ber \int(|\urre|^2+|\brre|^2)\le
C(\eps^{\kappa-1}+\eps_1^2+\eps_2^2+\frac{\eps_1+\eps_2}{\sqrt\eps}
+\frac{(\eps_1-\eps_2)^2}{\eps\sqrt\eps(\eps_1+\eps_2)}).\label{est-se5}\er
Now, combining \eqref{est-se1} and \eqref{est-se5}, we can get the
estimate \eqref{res1} in Theorem \ref{th1}.

For $L^\infty$ convergence rate, we just consider the case of
two-dimensional MHD system. For the three dimensional case, the
regularity problem involved here is open. To complete our Theorem,
we need to do higher order energy estimates. Of course, even though
the basic ideas of doing the higher order energy estimates are the
same as in $L^2$ estimates, but this is very complex. In the
following, we just consider the case $u^0=b^0$ and $0<z<\infty$. The
others are similar. First, to establish the convergence rate of high
order derivatives, we need solve the exact Prandtl's type equations
so as to obtain the more better convergence rate on the error
functions. Second, limited to the length of paper, we just give
estimates of some key terms which will appear when differentiating
nonlinear terms of MHD system and which are required to obtain the
estimates by using some kinds of different techniques from the
previous steps.

Let the boundary functions $u^\eps_B(x, z, t)=(u_{1B}^\eps,
u_{3B}^\eps)$ and $b^\eps_B(x, z, t)=(b_{1B}^\eps, b_{3B}^\eps)$
satisfy respectively the following Prandtl's type equations
\ber &&\pt u_{1B}^\eps=\eps\partial_z^2 u_{1B}^\eps\nonumber\\
&&\partial_1u^\eps_{1B}+\partial_3u^\eps_{3B}=0\nonumber\\
&& u_{1B}^\eps(x, z=0,t)=-u^0_1(x, z=0), u_{1B}^\eps(x, z=\infty,t)=0\nonumber\er and \ber &&\pt b_{1B}^\eps=\eps\partial_z^2 b_{1B}^\eps\nonumber\\
&&\partial_1b^\eps_{1B}+\partial_3b^\eps_{3B}=0\nonumber\\
&& b_{1B}^\eps(x, z=0,t)=-b^0_1(x, z=0), b_{1B}^\eps(x,
z=\infty,t)=0,\nonumber\er which can be solved as
\ber u_{1B}^\eps(x, z, t)=-u_1^0(x, z=0)\int_{\frac{z}{\sqrt{\eps (t+s)}}}^\infty \frac{e^{-{\xi^2}/4}}{\sqrt\pi}d\xi\nonumber\\
b_{1B}^\eps(x, z, t)=-b_1^0(x, z=0)\int_{\frac{z}{\sqrt{\eps
(t+s)}}}^\infty \frac{e^{-{\xi^2}/4}}{\sqrt\pi}d\xi\nonumber\er for
any given constant $s>0$ independent of $\eps$. It is easy to verify
that, due to the fact that $u_1^0(x, z=0)=b^0_1(x, z=0)=const$, \ber
\partial_tu^\eps_{3B}-\eps\partial_z^2u^\eps_{3B}=0,\partial_tb^\eps_{3B}-\eps\partial_z^2b^\eps_{3B}=0\nonumber\er
and the boundary functions $u^\eps_B, b^\eps_B$ have the similar
properties as given in Lemma \ref{lem1}. When $u^0=b^0$,
$u_B^\eps=b^\eps_B$.

Now, replacing $(\ue, \be)$ by $(u^{0}+\ubee+\urre,
b^{0}+\bbee+\brre)$ in the MHD system \eqref{vmhde1}-\eqref{vmhde5}
and using the system \eqref{imhd1}-\eqref{imhd5} in the
two-dimensional case, we gets that \ber
&&\pt\urre+\ue\cdot\nabla\urre+\urre\cdot\nabla\ubee+\urre\cdot\nabla
u^{0}
-\eps_1\Delta\urre-\eps_1\Delta u^{0}-(\eps_1-\eps)\partial_{z}^2\ubee\nonumber\\
&=&-\nabla
(\ppe-p^{0})+\be\cdot\nabla\brre+\brre\cdot\nabla\bbee+\brre\cdot\nabla  b^{0}, \quad \mbox{ in }\quad \Omega\times(0, T),\label{mhdre1c}\\
&&\pt\brre+\ue\cdot\nabla\brre+\urre\cdot\nabla\bbee+\urre\cdot\nabla b^{0}
-\eps_2\Delta\brre-\eps_2\Delta b^{0}-(\eps_2-\eps)\partial_{z}^2\bbee\nonumber\\
&=&\be\cdot\nabla\urre+\brre\cdot\nabla\ubee+\brre\cdot\nabla  b^{0}, \quad \mbox{ in }\quad \Omega\times(0, T),\label{mhdre2c}\\
&&div\ue=div\be=div\urre=div\brre=0,\mbox{ in } \Omega\times(0, T),\label{mhdre3c}\\
&&\urre=\brre=0,  \quad\mbox { on } \quad x\in \omega, z=0, 0\le
t\le
T\label{mhdre4c}\\
&&\urre(t=0)=\ue(0)-u^{0}(0)-\ubee(t=0),\nonumber\\
&&\brre(t=0)=\be(0)-  b^{0}(0)-\bbee(t=0), \quad \mbox{ on } \quad
\Omega .\label{mhdre5c}\er As before, \eqref{mhdre1c} and
\eqref{mhdre2c}  imply that \ber
&&\pt(\urre-\brre)+\ue\cdot\nabla(\urre-\brre)-\frac{\eps_1+\eps_2}2\Delta
(\urre-\brre)-\frac{\eps_1-\eps_2}2\Delta (\urre+\brre)
\nonumber\\
&=&-\nabla
(\ppe-p^{0})-\be\cdot\nabla(\urre-\brre)+\frac{\eps_1-\eps_2}2\Delta(u^\eps_B+b_B^\eps)+\frac{\eps_1-\eps_2}2\Delta(u^0+b^0)\label{mhdre-s1c}\er
Now, using \eqref{mhdre3c}-\eqref{mhdre4c} and taking the scalar
product of \eqref{mhdre-s1c} with $\urre-\brre$, yield that \ber
&&\frac 12\frac d{dt}\int
|\urre-\brre|^2+\frac{\eps_1+\eps_2}2\int|\nabla(\urre-\brre)|^2\nonumber\\
&=&-\frac{\eps_1-\eps_2}2\int\nabla((u_R^\eps+b_R^\eps)+(u_B^\eps+b_B^\eps)+(u^0+b^0))\cdot\nabla(u_R^\eps-b_R^\eps),\label{es-s1c}\er
which imply, by using \eqref{mhdre5c} and the properties of the
boundary layer functions and with the help of the Young's inequality
and the assumption \eqref{initial-e-assc} on initial data, that, for
$0\le t\le T<\infty$, \ber &&\int
|(\urre-\brre)(t)|^2+(1-\delta)(\eps_1+\eps_2)\int_0^t\int|\nabla(u_R^\eps-b_R^\eps)|^2\nonumber\\
&\le&
\int
|(\urre-\brre)(t=0)|^2+\frac{(\eps_1-\eps_2)^2}{4\delta(\eps_1+\eps_2)}[\int_0^t\int(|\nabla
u_R^\eps|^2+|\nabla b_R^\eps|^2)+Ct+\frac C{\sqrt\eps}t]\nonumber\\
&\le&
C\eps^\kappa+\frac{(\eps_1-\eps_2)^2}{4\delta(\eps_1+\eps_2)}\int_0^t\int(|\nabla
u_R^\eps|^2+|\nabla
b_R^\eps|^2)+C\frac{(\eps_1-\eps_2)^2}{4\delta\sqrt\eps(\eps_1+\eps_2)}\label{est-se1c}\er
for some constant $C>0$ and $\delta>0$ independent of
$\eps,\eps_1,\eps_2$, and $\kappa>2$ .

Now, by taking the scalar product of \eqref{mhdre1c} with $\urre$
and the scalar product of \eqref{mhdre2c} with $\brre$, we arrive at
\ber \frac 12\frac d{dt}\int
(|\urre|^2+|\brre|^2)+\eps_1\int|\nabla\urre|^2+\eps_2\int|\nabla\brre|^2
=\sum_{i=1}^{7}K_i,\label{est-se3c} \er where $K_i, i=1, \cdots, 7,$
are given respectively as follows \ber
&&K_1=-\int\nabla(\ppe-p^0)\urre;\nonumber\\
&&K_2=-\int\ue\cdot\nabla\urre\urre-\int\ue\cdot\nabla\brre\brre;\nonumber\\
&&K_3=-\int\urre\cdot\nabla\ubee\urre-\int\urre\cdot\nabla\bbee\brre +\int\brre\cdot\nabla\bbee\urre+\int\brre\cdot\nabla\ubee\brre\nonumber\\
&&K_4=-\int\urre\cdot\nabla u^{0}\urre-\int\urre\cdot\nabla  b^{0}\brre +\int\brre\cdot\nabla  b^{0}\urre+\int\brre\cdot\nabla u^{0}\brre; \nonumber\\
&&K_5=\int\be\cdot\nabla\brre\urre+\int\be\cdot\nabla\urre\brre;\nonumber\\
&&K_6=\eps_1\int\partial_z^2u^{0}\urre+\eps_2\int\partial_z^2
b^{0}\brre; \quad
K_7=(\eps_1-\eps)\int\partial_z^2\ubee\urre+(\eps_2-\eps)\int\partial_z^2
\bbee\brre.\nonumber\er We now bound each of $K_i, i=1, \cdots, 7$.

1) First,  $K_i, i=1,\cdots, 6$ can be estimated as before for $J_1,
J_3, J_6, J_8, J_9$.

2) Second, we compute \ber
K_7&=&-(\eps_1-\eps)\int\partial_z\ubee\partial_z\urre-(\eps_2-\eps)\int\partial_z
\bbee\partial_z\brre\nonumber\\
&\le& \delta\eps_1\int|\nabla u_R^\eps|^2+\delta\eps_2\int|\nabla
b_R^\eps|^2+C(\frac{(\eps_1-\eps)^2}{\eps_1\sqrt\eps}+\frac{(\eps_2-\eps)^2}{\eps_2\sqrt\eps}),\label{est-e10c}
\er  Combining \eqref{est-se3c} with \eqref{est-e1}-\eqref{est-e5}
for $K_i, i=1, \cdots, 6$ and \eqref{est-e10c} together and using
the assumptions \eqref{eps-assu} in Theorem \ref{th1}, we have \ber
&&\frac 12 \frac d{dt}\int(|\urre|^2+|\brre|^2)+(1-\delta)\eps_1\int
|\nabla\urre|^2+(1-\delta)\eps_2\int |\nabla\brre|^2\nonumber\\
&\le& C\int
(|\urre|^2+|\brre|^2)+\frac{(\eps_1-\eps_2)^2}{4\delta\eps(\eps_1+\eps_2)}\int_0^t\int(|\nabla
u_R^\eps|^2+|\nabla
b_R^\eps|^2)\nonumber\\
&&+C(\eps^{\kappa-1}+\eps_1^2+\eps_2^2+\frac{(\eps_1-\eps)^2}{\eps_1\sqrt\eps}+\frac{(\eps_2-\eps)^2}{\eps_2\sqrt\eps})
+C\frac{(\eps_1-\eps_2)^2}{4\delta\eps\sqrt\eps(\eps_1+\eps_2)}\nonumber\\
&\le& C[\int (|\urre|^2+|\brre|^2)+\eps_1\int_0^t\int|\nabla
u_R^\eps|^2+\eps_2\int_0^t\int|\nabla
b_R^\eps|^2]\nonumber\\
&&+C(\eps^{\kappa-1}+\eps_1^2+\eps_2^2+\frac{(\eps_1-\eps)^2}{\eps_1\sqrt\eps}+\frac{(\eps_2-\eps)^2}{\eps_2\sqrt\eps})
+C\frac{(\eps_1-\eps_2)^2}{4\delta\eps\sqrt\eps(\eps_1+\eps_2)}.\label{est-se4c}\er
It follows from \eqref{est-se4c} and by Gronwall's inequality that
\ber \int(|\urre|^2+|\brre|^2)+\eps_1\int_0^t\int|\nabla
u_R^\eps|^2+\eps_2\int_0^t\int|\nabla b_R^\eps|^2 \le
C\beta_0(\eps).\label{est-se5c}\er Here \ber
\beta_0(\eps)=\eps^{\kappa-1}+\eps_1^2+\eps_2^2+\frac{(\eps_1-\eps)^2}{\eps_1\sqrt\eps}+\frac{(\eps_2-\eps)^2}{\eps_2\sqrt\eps}
+\frac{(\eps_1-\eps_2)^2}{\eps\sqrt\eps(\eps_1+\eps_2)}.\label{beta0}\er
Then we have \ber &&\int
|(\urre-\brre)(t)|^2+(1-\delta)(\eps_1+\eps_2)\int_0^t\int|\nabla(u_R^\eps-b_R^\eps)|^2\nonumber\\
&\le&
C\frac{(\eps_1-\eps_2)^2}{(\eps_1+\eps_2)\min\{\eps_1,\eps_2\}}\beta_0(\eps)
+C\eps^\kappa+C\frac{(\eps_1-\eps_2)^2}{\sqrt\eps(\eps_1+\eps_2)}=\bar\beta_0(\eps)\label{est-se1c-1}\er
Differentiate the equation \eqref{mhdre-s1c} in time, multiply the
resulting one by $\partial_t(u_R^\eps-b^\eps_R)$ and integrate over
$\Omega$. Notice that $\partial_tu^\eps_R|_{t=0}=\eps_1\Delta
u^0(t=0)+O(\eps^\kappa+\eps_1\eps^\kappa+\frac{\eps_1-\eps}\eps)$,
$\partial_tb^\eps_R|_{t=0}=\eps_2\Delta
b^0(t=0)+O(\eps^\kappa+\eps_2\eps^\kappa+\frac{\eps_2-\eps}\eps)$
and
\ber &&|\int \partial_t(u^\eps+b^\eps)\cdot \nabla(u_R^\eps-b^\eps_R)\cdot\partial_t(u_R^\eps-b^\eps_R)|\nonumber\\
&\le& |\int \partial_t(u^\eps_R+b^\eps_R)\cdot \nabla(u_R^\eps-b^\eps_R)\cdot\partial_t(u_R^\eps-b^\eps_R)|\nonumber\\
&&+|\int \partial_t(u^0+b^0+u^\eps_B+b_B^\eps)\cdot \nabla(u_R^\eps-b^\eps_R)\cdot\partial_t(u_R^\eps-b^\eps_R)|\nonumber\\
&\le&\|\nabla(u^\eps_R-b^\eps_R)\|_{L^2}\|\partial_t(u^\eps_R+b^\eps_R)\partial_t(u_R^\eps-b^\eps_R)\|_{L^2}
+C\|\nabla(u^\eps_R-b^\eps_R)\|_{L^2}\|\partial_t(u_R^\eps-b^\eps_R)\|_{L^2}\nonumber\\
&\le&\|\nabla(u^\eps_R-b^\eps_R)\|_{L^2}(\|\partial_t(u^\eps_R+b^\eps_R)\|_{L^4}^2+\|\partial_t(u_R^\eps-b^\eps_R)\|_{L^4}^2)\nonumber\\
&&+C\|\nabla(u^\eps_R-b^\eps_R)\|_{L^2}\|\partial_t(u_R^\eps-b^\eps_R)\|_{L^2}\nonumber\\
&\le&\|\nabla(u^\eps_R-b^\eps_R)\|_{L^2}(\|\partial_t(u^\eps_R+b^\eps_R)\|_{L^2}\|\nabla\partial_t(u^\eps_R+b^\eps_R)\|_{L^2}
+\|\partial_t(u^\eps_R-b^\eps_R)\|_{L^2}\|\nabla\partial_t(u^\eps_R-b^\eps_R)\|_{L^2})\nonumber\\
&&+C\|\nabla(u^\eps_R-b^\eps_R)\|_{L^2}^2+C\|\partial_t(u_R^\eps-b^\eps_R)\|_{L^2}^2\nonumber\\
&\le&\delta(\eps_1+\eps_2)\|\nabla\partial_t(u_R^\eps-b^\eps_R)\|_{L^2}^2+\delta\eps^2\|\partial_t\nabla (u^\eps_R+b^\eps_R)\|_{L^2}^2\nonumber\\
&&+C\frac1{\eps_1+\eps_2}\|\nabla(u^\eps_R-b^\eps_R)\|_{L^2}^2\|\partial_t(u^\eps_R-b^\eps_R)\|_{L^2}^2
+C\frac1{\eps^2}\|\nabla(u^\eps_R-b^\eps_R)\|_{L^2}^2\|\partial_t(u^\eps_R+b^\eps_R)\|_{L^2}^2\nonumber\\
&&+C\|\nabla(u^\eps_R-b^\eps_R)\|_{L^2}^2+C\|\partial_t(u_R^\eps-b^\eps_R)\|_{L^2}^2\nonumber\er
Hence we have
 \ber &&\frac d{dt}\int
|\partial_t(\urre-\brre)(t)|^2+(1-\delta)(\eps_1+\eps_2)\int|\partial_t\nabla(u_R^\eps-b_R^\eps)|^2\nonumber\\
&\le&C\frac1{\eps_1+\eps_2}\|\nabla(u^\eps_R-b^\eps_R)\|_{L^2}^2\|\partial_t(u^\eps_R-b^\eps_R)\|_{L^2}^2\nonumber\\
&&+C\frac1{\eps^2}\|\nabla(u^\eps_R-b^\eps_R)\|_{L^2}^2\|\partial_t(u^\eps_R+b^\eps_R)\|_{L^2}^2\nonumber\\
&&+C\|\nabla(u^\eps_R-b^\eps_R)\|_{L^2}^2+C\|\partial_t(u_R^\eps-b^\eps_R)\|_{L^2}^2\nonumber\\
&&+C(\frac{(\eps_1-\eps_2)^2}{4\delta(\eps_1+\eps_2)}+\delta\eps^2)\int(|\partial_t\nabla
u_R^\eps|^2+|\partial_t\nabla
b_R^\eps|^2)+C\frac{(\eps_1-\eps_2)^2}{4\delta\sqrt\eps(\eps_1+\eps_2)}\label{est-se1ct}\er
or \ber &&\int
|\partial_t(\urre-\brre)(t)|^2+(1-\delta)(\eps_1+\eps_2)\int_0^t\int|\partial_t\nabla(u_R^\eps-b_R^\eps)|^2\nonumber\\
&\le&C\int_0^t\frac{\|\nabla(u^\eps_R-b^\eps_R)\|_{L^2}^2}{\eps^2}\|\partial_t(u^\eps_R+b^\eps_R)\|_{L^2}^2+
C\frac{(\eps_1-\eps_2)^2}{(\eps_1+\eps_2)^2\min\{\eps_1,\eps_2\}}\beta_0(\eps)\nonumber\\
&&+C\frac{\eps^\kappa}{\eps_1+\eps_2}+C\frac{(\eps_1-\eps_2)^2}{\sqrt\eps(\eps_1+\eps_2)^2}+C(\eps_1-\eps_2)^2\nonumber\\
&&+(\frac{C(\eps_1-\eps_2)^2}{4\delta(\eps_1+\eps_2)}+\delta\eps^2)\int_0^t\int(|\partial_t\nabla
u_R^\eps|^2+|\partial_t\nabla
b_R^\eps|^2)+C\frac{(\eps_1-\eps_2)^2}{4\delta\sqrt\eps(\eps_1+\eps_2)}\label{est-se1ct-1}\er
for some $\delta>0$ independent of $\eps$. Here we require
$\frac{\beta_0(\eps)}{\min\{\eps_1,\eps_2\}(\eps_1+\eps_2)}\le C$.

Differentiate the equations \eqref{mhdre1c} and \eqref{mhdre2c} in
time, multiply the resulting ones respectively by
$\partial_tu_R^\eps$ and $\partial_t b^\eps_R$ and integrate over
$\Omega$. Notice that \ber &&|\int (-\partial_tu^\eps\cdot\nabla
u^\eps_R\partial_tu^\eps_R+\partial_tb^\eps\cdot\nabla
b^\eps_R\partial_tu^\eps_R-\partial_tu^\eps\cdot\nabla
b^\eps_R\partial_tb^\eps_R+\partial_tb^\eps\cdot\nabla
u^\eps_R\partial_tb^\eps_R)|\nonumber\\
&\le&|\int (-\partial_tu^\eps_R\cdot\nabla
u^\eps_R\partial_tu^\eps_R+\partial_tb^\eps_R\cdot\nabla
b^\eps_R\partial_tu^\eps_R-\partial_tu^\eps_R\cdot\nabla
b^\eps_R\partial_tb^\eps_R+\partial_tb^\eps_R\cdot\nabla
u^\eps_R\partial_tb^\eps_R)|\nonumber\\
&&+|\int (-\partial_t(u^0+u_B^\eps)\cdot\nabla
u^\eps_R\partial_tu^\eps_R+\partial_t(b^0+b_B^\eps)\cdot\nabla
b^\eps_R\partial_tu^\eps_R\nonumber\\
&&\quad-\partial_t(u^0+u_B^\eps)\cdot\nabla
b^\eps_R\partial_tb^\eps_R+\partial_t(b^0+b^\eps_B)\cdot\nabla
u^\eps_R\partial_tb^\eps_R)|\nonumber\\
&\le&C\|\nabla
u^\eps_R\|_{L^2}(\|\partial_tu^\eps_R\|_{L^4}^2+\|\partial_tb^\eps_R\|_{L^4}^2)+C\|\nabla
b^\eps_R\|_{L^2}\|\partial_tu^\eps_R\|_{L^4}^2+\|\partial_tb^\eps_R\|_{L^4}^2\nonumber\\
&&\quad+C(\|\nabla u^\eps_R\|^2_{L^2}+\|\nabla
b^\eps_R\|^2_{L^2})+C\|\partial_tu^\eps_R\|_{L^2}^2+C\|\partial_tb^\eps_R\|_{L^2}^2\nonumber\\
&\le&C(\|\nabla u^\eps_R\|_{L^2}+\|\nabla
b^\eps_R\|_{L^2})(\|\partial_tu^\eps_R\|_{L^2}\|\nabla\partial_tu^\eps_R\|_{L^2}+\|\partial_tb^\eps_R\|_{L^2}\|\nabla\partial_tb^\eps_R\|_{L^2})\nonumber\\
&&\quad+C(\|\nabla u^\eps_R\|^2_{L^2}+\|\nabla
b^\eps_R\|^2_{L^2})+C\|\partial_tu^\eps_R\|_{L^2}^2+C\|\partial_tb^\eps_R\|_{L^2}^2\nonumber\\
&\le&\delta\eps_1\|\nabla\partial_tu^\eps_R\|_{L^2}^2+C\frac{\|\nabla
u^\eps_R\|_{L^2}^2}{\eps_1}\|\partial_tu^\eps_R\|_{L^2}^2+\delta\eps_2\|\nabla\partial_tb^\eps_R\|_{L^2}^2+C\frac{\|\nabla
u^\eps_R\|_{L^2}^2}{\eps_2}\|\partial_tb^\eps_R\|_{L^2}^2\nonumber\\
&&+C\frac{\|\nabla
b^\eps_R\|_{L^2}^2}{\eps_2}\|\partial_tb^\eps_R\|_{L^2}^2+C\frac{\|\nabla
b^\eps_R\|_{L^2}^2}{\eps_1}\|\partial_tu^\eps_R\|_{L^2}^2\nonumber\\
&&+C(\|\nabla u^\eps_R\|^2_{L^2}+\|\nabla
b^\eps_R\|^2_{L^2}+\|\partial_tu^\eps_R\|_{L^2}^2+\|\partial_tb^\eps_R\|_{L^2}^2).\label{est-pt1}\er
\ber &&|\int(-\partial_tu^\eps_R\cdot\nabla u^\eps_B\partial_t
u^\eps_R+\partial_tb^\eps_R\cdot\nabla
b^\eps_B\partial_tu^\eps_R-\partial_tu^\eps_R\nabla
b^\eps_B\partial_tb^\eps_R+\partial_tb^\eps_R\cdot\nabla
u^\eps_B\partial_tb^\eps_R|\nonumber\\
&=&|-\int\partial_t(u^\eps_R-b^\eps_R)\cdot\nabla
u^\eps_B\cdot\partial_t(u^\eps_R+b^\eps_R)|\nonumber\\
&\le&C\int(|\partial_tu^\eps_R|^2+|\partial_tu^\eps_R|^2)+C\int\frac{|\partial_t(u^\eps_R-b^\eps_R)|^2}{\eps}.\label{est-pt2}\er
\ber &&|\int(-u^\eps_R\cdot\partial_t\nabla u^\eps_B\partial_t
u^\eps_R+b^\eps_R\cdot\nabla
\partial_tb^\eps_B\partial_tu^\eps_R-u^\eps_R\nabla\partial_t
b^\eps_B\partial_tb^\eps_R+b^\eps_R\cdot\nabla
\partial_tu^\eps_B\partial_tb^\eps_R|\nonumber\\
&=&|-\int(u^\eps_R-b^\eps_R)\cdot\nabla\partial_t
u^\eps_B\cdot\partial_t(u^\eps_R+b^\eps_R)|\nonumber\\
&\le&C\int(|\partial_tu^\eps_R|^2+|\partial_tu^\eps_R|^2)+C\int\frac{|u^\eps_R-b^\eps_R|^2}{\eps}.\label{est-pt2-1}\er
Hence we have \ber &&\frac d{dt}\int(|\partial_t
u^\eps_R|^2+|\partial_tb^\eps_R|^2)+\eps_1\int|\nabla\partial_tu^\eps_R|^2+\eps_2\int|\nabla\partial_tb^\eps_R|^2\nonumber\\
&\le&C\frac{\|\nabla
u^\eps_R\|_{L^2}^2}{\eps_1}\int|\partial_tu^\eps_R|^2+C\frac{\|\nabla
u^\eps_R\|_{L^2}^2}{\eps_2}\int|\partial_tb^\eps_R|^2\nonumber\\
&&+C\frac{\|\nabla
b^\eps_R\|_{L^2}^2}{\eps_2}\int|\partial_tb^\eps_R|^2+C\frac{\|\nabla
b^\eps_R\|_{L^2}^2}{\eps_1}\int|\partial_tu^\eps_R|^2\nonumber\\
&&+C\int(|\nabla u^\eps_R|^2+|\nabla
b^\eps_R|^2)+C\int(|\partial_tu^\eps_R|^2+|\partial_tb^\eps_R|^2)\nonumber\\
&&+C\int(|u^\eps_R|^2+|b^\eps_R|^2)+C\int\frac{|\partial_t(u^\eps_R-b^\eps_R)|^2+|(u^\eps_R-b^\eps_R)|^2}{\eps}\nonumber\\
&&+C(\frac{(\eps_1-\eps)^2}{\eps_1\sqrt\eps}+\frac{(\eps_2-\eps)^2}{\eps_1\sqrt\eps})\nonumber\\
&\le&C\frac{\|\nabla
u^\eps_R\|_{L^2}^2}{\eps_1}\int|\partial_tu^\eps_R|^2+C\frac{\|\nabla
u^\eps_R\|_{L^2}^2}{\eps_2}\int|\partial_tb^\eps_R|^2\nonumber\\
&&+C\frac{\|\nabla
b^\eps_R\|_{L^2}^2}{\eps_2}\int|\partial_tb^\eps_R|^2+C\frac{\|\nabla
b^\eps_R\|_{L^2}^2}{\eps_1}\int|\partial_tu^\eps_R|^2\nonumber\\
&&+C\int(|\nabla u^\eps_R|^2+|\nabla
b^\eps_R|^2)+C\int(|\partial_tu^\eps_R|^2+|\partial_tb^\eps_R|^2)\nonumber\\
&&+C\int(|u^\eps_R|^2+|b^\eps_R|^2)
+C\int_0^t\frac{\|\nabla(u^\eps_R-b^\eps_R)\|_{L^2}^2}{\eps^3}\int (|\partial_tu^\eps_R|^2+|\partial_tb^\eps_R|^2)\nonumber\\
&&+C(\frac{(\eps_1-\eps_2)^2}{\eps(\eps_1+\eps_2)}+\delta\eps)\int_0^t\int(|\partial_t\nabla
u_R^\eps|^2+|\partial_t\nabla
b_R^\eps|^2)+C\frac{(\eps_1-\eps_2)^2}\eps+C\eps^{\kappa-1}\nonumber\\
&&+C\frac{(\eps_1-\eps_2)^2}{\eps\sqrt\eps(\eps_1+\eps_2)}
+C\frac{(\eps_1-\eps_2)^2}{\eps(\eps_1+\eps_2)\min\{\eps_1,\eps_2\}}\beta_0(\eps)+C(\frac{(\eps_1-\eps)^2}{\eps_1\sqrt\eps}
+\frac{(\eps_2-\eps)^2}{\eps_1\sqrt\eps})\nonumber\\
&&+C\frac{(\eps_1-\eps_2)^2}{\eps(\eps_1+\eps_2)^2\min\{\eps_1,\eps_2\}}\beta_0(\eps)+C\frac{\eps^{\kappa-1}}{\eps_1+\eps_2}\nonumber\\
&&+C\frac{(\eps_1-\eps_2)^2}{\eps\sqrt\eps(\eps_1+\eps_2)^2}+C(\eps_1-\eps_2)^2.\label{est-pt3}
\er It follows from Gronwall's inequality and \eqref{est-pt3}, with
the help of the estimates \eqref{est-se5c} and \eqref{est-se1c-1},
that \ber\int(|\partial_t
u^\eps_R|^2+|\partial_tb^\eps_R|^2)+\eps_1\int_0^t\int|\nabla\partial_tu^\eps_R|^2+\eps_2\int_0^t\int|\nabla\partial_tb^\eps_R|^2\le
C\beta_1(\eps).\label{est-pt4}\er Here \ber
\beta_1(\eps)=\eps_1^2+\eps^2_2+\frac{(\eps_1-\eps_2)^2}{\eps(\eps_1+\eps_2)^2\min\{\eps_1,\eps_2\}}\beta_0(\eps)+\frac{\eps^{\kappa-1}}{\eps_1+\eps_2}\nonumber\\
+\frac{(\eps_1-\eps_2)^2}{\eps\sqrt\eps(\eps_1+\eps_2)^2}+(\frac{(\eps_1-\eps)^2}{\eps_1\sqrt\eps}+\frac{(\eps_2-\eps)^2}{\eps_1\sqrt\eps}).\label{est-pt-beta-1}\er
Here we require $\frac{\bar\beta_0(\eps)}{\eps_1+\eps_2)\eps^3}\le
C.$

Then it follows from \eqref{est-se1ct} and \eqref{est-pt4} that
\ber\int
|\partial_t(\urre-\brre)(t)|^2+(1-\delta)(\eps_1+\eps_2)\int_0^t\int|\partial_t\nabla(u_R^\eps-b_R^\eps)|^2\le
C\bar\beta_1(\eps),\nonumber\er
\ber\bar\beta_1(\eps)=(\eps_1-\eps_2)^2+\frac{(\eps_1-\eps_2)^2}{(\eps_1+\eps_2)\min\{\eps_1,\eps_2\}}\beta_1(\eps)
+\frac{\eps^2}{\min\{\eps_1,\eps_2\}}\beta_1(\eps)+\frac{(\eps_1-\eps_2)^2}{\sqrt\eps(\eps_1+\eps_2)}.\label{est-se1ct-1}\er
Similarly, we can obtain the estimates on tangential derivatives as
follows:\ber\int(|\partial_x
u^\eps_R|^2+|\partial_xb^\eps_R|^2)+\eps_1\int_0^t\int|\nabla\partial_xu^\eps_R|^2+\eps_2\int_0^t\int|\nabla\partial_xb^\eps_R|^2\le
C\beta_2(\eps).\label{est-pt5}\er Here \ber
\beta_2(\eps)=\frac{(\eps_1-\eps_2)^2}{\eps(\eps_1+\eps_2)^2\min\{\eps_1,\eps_2\}}\beta_0(\eps)+\eps_1^2+\eps^2_2
+\frac{\eps^{\kappa-1}}{\eps_1+\eps_2}+\frac{(\eps_1-\eps_2)^2}{\eps\sqrt\eps(\eps_1+\eps_2)^2}.\label{est-pt-beta-2}\er
and \ber \int
|\partial_x(\urre-\brre)(t)|^2+(\eps_1+\eps_2)\int_0^t\int|\nabla\partial_x(u_R^\eps-b_R^\eps)|^2\le
C\bar\beta_2(\eps),\er
\ber\bar\beta_2(\eps)=\eps^\kappa+\frac{(\eps_1-\eps_2)^2}{(\eps_1+\eps_2)\min\{\eps_1,\eps_2\}}\beta_2(\eps)
+\frac{\eps^2}{\min\{\eps_1,\eps_2\}}\beta_2(\eps)+\frac{(\eps_1-\eps_2)^2}{\sqrt\eps(\eps_1+\eps_2)}.\label{est-se1ct-2}\er
Here $\|\partial_x(u^\eps_R(t=0),
b^\eps_R(t=0))\|_{L^2}^2=O(\eps^\kappa)$ and
$\partial_xu^\eps_B=\partial_xb^\eps_B=0$.

Finally, we apply $\partial_t\partial_x$ to the equations
\eqref{mhdre1c}-\eqref{mhdre4c} and \eqref{mhdre-s1c}, multiply the
resulting ones respectively by
$\partial_t\partial_xu^\eps_R,\partial_t\partial_xb^\eps_R$ and
$\partial_t\partial_x(u^\eps_R-b^\eps_R)$, and integrate over
$\Omega$. Notice that $\|\partial_t\partial_x(u^\eps_R(t=0),
b^\eps_R(t=0))\|_{L^2}^2=O(\eps^\kappa)$,
$\eps_1\partial_t\partial_x\Delta
u^0=0=\eps_2\partial_t\partial_x\Delta b^0$ and
$\partial_t\partial_xu^\eps_B=\partial_t\partial_xb^\eps_B=0$. Also,
\ber &&|\int \partial_t\partial_x(u^\eps+b^\eps)\cdot \nabla(u_R^\eps-b^\eps_R)\cdot\partial_t\partial_x(u_R^\eps-b^\eps_R)|\nonumber\\
&=& |\int \partial_t\partial_x(u^\eps_R+b^\eps_R)\cdot \nabla(u_R^\eps-b^\eps_R)\cdot\partial_t\partial_x(u_R^\eps-b^\eps_R)|\nonumber\\
&\le&\|\nabla(u^\eps_R-b^\eps_R)\|_{L^2}\|\partial_t\partial_x(u^\eps_R+b^\eps_R)\partial_t\partial_x(u_R^\eps-b^\eps_R)\|_{L^2}\nonumber\\
&\le&C\|\nabla(u^\eps_R-b^\eps_R)\|_{L^2}(\|\partial_t\partial_x(u^\eps_R+b^\eps_R)\|_{L^4}^2+\|\partial_t\partial_x(u_R^\eps-b^\eps_R)\|_{L^4}^2)\nonumber\\
&\le&C\|\nabla(u^\eps_R-b^\eps_R)\|_{L^2}(\|\partial_t\partial_x(u^\eps_R+b^\eps_R)\|_{L^2}^2\|\nabla\partial_t\partial_x(u^\eps_R+b^\eps_R)\|_{L^2}^2\nonumber\\
&&+\|\partial_t\partial_x(u_R^\eps-b^\eps_R)\|_{L^2}^2\|\nabla\partial_t\partial_x(u_R^\eps-b^\eps_R)\|_{L^2}^2)\nonumber\\
&\le&\delta(\eps_1+\eps_2)\|\nabla\partial_t\partial_x(u_R^\eps-b^\eps_R)\|_{L^2}^2+\delta\eps^2\|\nabla\partial_t\partial_x(u_R^\eps+b^\eps_R)\|_{L^2}^2\nonumber\\
&&+C\frac{\|\nabla(u^\eps_R-b^\eps_R)\|_{L^2}^2}{\eps_1+\eps_2}\|\partial_t\partial_x(u^\eps_R-b^\eps_R)\|_{L^2}^2
+C\frac{\|\nabla(u^\eps_R-b^\eps_R)\|_{L^2}^2}{\eps^2}\|\nabla\partial_t\partial_x(u^\eps_R+b^\eps_R)\|_{L^2}^2\label{est-pt6}\er
and
\ber &&|\int \partial_t(u^\eps+b^\eps)\cdot \nabla\partial_x(u_R^\eps-b^\eps_R)\cdot\partial_t\partial_x(u_R^\eps-b^\eps_R)|\nonumber\\
&\le& |\int \partial_t(u^\eps_R+b^\eps_R)\cdot \nabla\partial_x(u_R^\eps-b^\eps_R)\cdot\partial_t\partial_x(u_R^\eps-b^\eps_R)|\nonumber\\
&&+|\int \partial_t(u^0+b^0+u^\eps_B+b_B^\eps)\cdot \nabla\partial_x(u_R^\eps-b^\eps_R)\cdot\partial_t\partial_x(u_R^\eps-b^\eps_R)|\nonumber\\
&\le&\|\nabla\partial_x(u^\eps_R-b^\eps_R)\|_{L^2}\|\partial_t(u^\eps_R+b^\eps_R)\partial_t\partial_x(u_R^\eps-b^\eps_R)\|_{L^2}
+C\|\nabla\partial_x(u^\eps_R-b^\eps_R)\|_{L^2}\|\partial_t\partial_x(u_R^\eps-b^\eps_R)\|_{L^2}\nonumber\\
&\le&C\|\nabla\partial_x(u^\eps_R-b^\eps_R)\|_{L^2}(\|\partial_t(u^\eps_R+b^\eps_R)\|_{L^4}^2+\|\partial_t\partial_x(u_R^\eps-b^\eps_R)\|_{L^4}^2)\nonumber\\
&&+C\|\nabla\partial_x(u^\eps_R-b^\eps_R)\|_{L^2}\|\partial_t\partial_x(u_R^\eps-b^\eps_R)\|_{L^2}\nonumber\\
&\le&C\|\nabla\partial_x(u^\eps_R-b^\eps_R)\|_{L^2}(\|\partial_t(u^\eps_R+b^\eps_R)\|_{L^2}^2\|\nabla\partial_t(u^\eps_R+b^\eps_R)\|_{L^2}^2
+\|\partial_t\partial_x(u_R^\eps-b^\eps_R)\|_{L^2}^2\|\nabla\partial_t\partial_x(u_R^\eps-b^\eps_R)\|_{L^2}^2)\nonumber\\
&&+C\|\nabla\partial_x(u^\eps_R-b^\eps_R)\|_{L^2}\|\partial_t\partial_x(u_R^\eps-b^\eps_R)\|_{L^2}\nonumber\\
&\le&\delta(\eps_1+\eps_2)\|\nabla\partial_t\partial_x(u_R^\eps-b^\eps_R)\|_{L^2}^2
+C\frac{\|\nabla\partial_x(u^\eps_R-b^\eps_R)\|_{L^2}^2}{\eps_1+\eps_2}\|\partial_t\partial_x(u^\eps_R-b^\eps_R)\|_{L^2}^2\nonumber\\
&&+\delta\eps^2\|\nabla\partial_t(u_R^\eps+b^\eps_R)\|_{L^2}^2
+C\frac{\|\nabla\partial_x(u^\eps_R-b^\eps_R)\|_{L^2}^2}{\eps^2}\|\partial_t(u^\eps_R+b^\eps_R)\|_{L^2}^2\nonumber\\
&&+C\|\nabla\partial_x(u^\eps_R-b^\eps_R)\|_{L^2}^2+C\|\partial_t\partial_x(u_R^\eps-b^\eps_R)\|_{L^2}^2.\label{est-pt7}\er
Similarly,$\int \partial_x(u^\eps+b^\eps)\cdot
\nabla\partial_t(u_R^\eps-b^\eps_R)\cdot\partial_t\partial_x(u_R^\eps-b^\eps_R)$can
be controlled by
$\delta(\eps_1+\eps_2)\|\nabla\partial_t\partial_x(u_R^\eps-b^\eps_R)\|_{L^2}^2+\delta\eps^2\|\nabla\partial_x(u_R^\eps+b^\eps_R)\|_{L^2}^2
+C\frac{\|\nabla\partial_t(u^\eps_R-b^\eps_R)\|_{L^2}^2}{\eps_1+\eps_2}\|\partial_t\partial_x(u^\eps_R-b^\eps_R)\|_{L^2}^2
+C\frac{\|\nabla\partial_t(u^\eps_R-b^\eps_R)\|_{L^2}^2}{\eps^2}\|\partial_x(u^\eps_R+b^\eps_R)\|_{L^2}^2
+C\|\nabla\partial_t(u^\eps_R-b^\eps_R)\|_{L^2}^2+C\|\partial_t\partial_x(u_R^\eps-b^\eps_R)\|_{L^2}^2$.
Thus, we have
\ber&&\frac d{dt}\int|\partial_t\partial_x(u^\eps_R-b_R^\eps)|^2+(1-\delta)(\eps_1+\eps_2)\int|\nabla\partial_t\partial_x(u^\eps_R-b^\eps_R)|^2\nonumber\\
&\le&C(\frac{(\eps_1-\eps_2)^2}{\eps_1+\eps_2}+\delta\eps^2)\int(|\nabla\partial_t\partial_xu^\eps_R|^2+|\nabla\partial_t\partial_xb^\eps_R|^2)\nonumber\\
&&+C\int|\nabla\partial_x(u^\eps_R-b^\eps_R)|^2+C\int|\nabla\partial_t(u^\eps_R-b^\eps_R)|^2+\int|\partial_t\partial_x(u^\eps_R-b^\eps_R)|^2\nonumber\\
&&+C\frac{\|\nabla\partial_t(u^\eps_R-b^\eps_R)\|_{L^2}^2}{\eps_1+\eps_2}\int |\partial_t\partial_x(u^\eps_R-b^\eps_R)|^2\nonumber\\
&&+C\frac{\|\nabla\partial_x(u^\eps_R-b^\eps_R)\|_{L^2}^2}{\eps_1+\eps_2}\int |\partial_t\partial_x(u^\eps_R-b^\eps_R)|^2\nonumber\\
&&+C\frac{\|\nabla\partial_x(u^\eps_R-b^\eps_R)\|_{L^2}^2}{\eps^2}\int |\partial_t(u^\eps_R+b^\eps_R)|^2\nonumber\\
&&+C\frac{\|\nabla\partial_t(u^\eps_R-b^\eps_R)\|_{L^2}^2}{\eps^2}\int
|\partial_x(u^\eps_R+b^\eps_R)|^2\label{est-pt8} \er or
\ber&&\int|\partial_t\partial_x(u^\eps_R-b_R^\eps)|^2+(1-\delta)(\eps_1+\eps_2)\int_0^t\int|\nabla\partial_t\partial_x(u^\eps_R-b^\eps_R)|^2\nonumber\\
&\le&C(\frac{(\eps_1-\eps_2)^2}{\eps_1+\eps_2}+\delta\eps^2)\int_0^t\int(|\nabla\partial_t\partial_xu^\eps_R|^2
+|\nabla\partial_t\partial_xb^\eps_R|^2)+C\beta_3(\eps).\label{est-pt8}
\er Here \ber&&\beta_3(\eps)=\eps^\kappa+\frac
{\beta_1(\eps)\bar\beta_2(\eps)+\beta_2(\eps)\bar\beta_1(\eps)}{\eps^2(\eps_1+\eps_2)}.\label{beta3}\er
Here we require $\frac{\bar\beta_1(\eps)}{(\eps_1+\eps_2)^2}\le C.$

On the other hand, \ber
&&|-\int(\partial_x(u^\eps_R-b^\eps_R)\cdot\nabla\partial_tu^\eps_B+\partial_t\partial_x(u^\eps_R-b^\eps_R)\cdot\nabla
u^\eps_B)\cdot\partial_t\partial_x(u^\eps_R+b^\eps_R)|\nonumber\\
&\le&C\int(|\partial_t\partial_xu^\eps_R|^2+|\partial_t\partial_xb^\eps_R|^2)
+C\int\frac{|\partial_x(u^\eps_R-b^\eps_R)|^2+|\partial_t\partial_x(u^\eps_R-b^\eps_R)|^2}{\eps}\nonumber\\
&\le&C\int(|\partial_t\partial_xu^\eps_R|^2+|\partial_t\partial_xb^\eps_R|^2)
+C(\frac{(\eps_1-\eps_2)^2}{\eps(\eps_1+\eps_2)}+\delta\eps)\int_0^t\int(|\nabla\partial_t\partial_xu^\eps_R|^2+|\nabla\partial_t\partial_xb^\eps_R|^2)\nonumber\\
&&+C\frac{\beta_3(\eps)}\eps
+C\eps^{\kappa-1}+C\frac{(\eps_1-\eps_2)^2}{\eps(\eps_1+\eps_2)\min\{\eps_1,\eps_2\}}\beta_2(\eps)
+C\frac{(\eps_1-\eps_2)^2}{\eps\sqrt\eps(\eps_1+\eps_2)}.\label{est-pt9}\er
As in \eqref{est-pt6} and \eqref{est-pt7}, one can get that \ber
&&\int[-(\partial_t\partial_xu^\eps\cdot\nabla
u^\eps_R+\partial_tu^\eps\cdot\nabla\partial_xu^\eps_R+\partial_xu^\eps\cdot\nabla\partial_tu^\eps_R)\nonumber\\
&&+(\partial_t\partial_xb^\eps\cdot\nabla
b^\eps_R+\partial_tb^\eps\cdot\nabla\partial_xb^\eps_R+\partial_xb^\eps\cdot\nabla\partial_tb^\eps_R)]\cdot\partial_t\partial_xu^\eps_R\nonumber\\
&\le&\delta\eps_1\|\nabla\partial_t\partial_xu^\eps_R\|^2_{L^2}+C\frac{\|\nabla
u^\eps_R\|^2}{\eps_1}\|\partial_t\partial_xu^\eps_R\|^2_{L^2}+C\frac{\|\nabla\partial_x
u^\eps_R\|^2}{\eps_1}\|\partial_tu^\eps_R\|^2_{L^2}\nonumber\\
&&+C\frac{\|\nabla\partial_t
u^\eps_R\|^2}{\eps_1}\|\partial_xu^\eps_R\|^2_{L^2}+C\|\nabla\partial_tu^\eps_R\|^2_{L^2}+C\|\nabla\partial_xu^\eps_R\|^2_{L^2}+C\frac{\|\nabla
b^\eps_R\|^2}{\eps_1}\|\partial_t\partial_xb^\eps_R\|^2_{L^2}\nonumber\\
&&+C\frac{\|\nabla\partial_x
b^\eps_R\|^2}{\eps_1}\|\partial_tb^\eps_R\|^2_{L^2}+C\frac{\|\nabla\partial_t
b^\eps_R\|^2}{\eps_1}\|\partial_xb^\eps_R\|^2_{L^2}\nonumber\\
&&+C\|\nabla\partial_tb^\eps_R\|^2_{L^2}+C\|\nabla\partial_xb^\eps_R\|^2_{L^2}+C\|\nabla
u^\eps_R\|^2_{L^2}+C\|\nabla
b^\eps_R\|^2_{L^2}.\label{est-pt10-1}\er \ber
&&\int[-(\partial_t\partial_xu^\eps\cdot\nabla
b^\eps_R+\partial_tu^\eps\cdot\nabla\partial_xb^\eps_R+\partial_xu^\eps\cdot\nabla\partial_tb^\eps_R)\nonumber\\
&&+(\partial_t\partial_xb^\eps\cdot\nabla
u^\eps_R+\partial_tb^\eps\cdot\nabla\partial_xu^\eps_R+\partial_xb^\eps\cdot\nabla\partial_tu^\eps_R)]\cdot\partial_t\partial_xb^\eps_R\nonumber\\
&\le&\delta\eps_2\|\nabla\partial_t\partial_xb^\eps_R\|^2_{L^2}+C\frac{\|\nabla
b^\eps_R\|^2}{\eps_2}\|\partial_t\partial_xu^\eps_R\|^2_{L^2}+C\frac{\|\nabla\partial_x
b^\eps_R\|^2}{\eps_2}\|\partial_tu^\eps_R\|^2_{L^2}\nonumber\\
&&+C\frac{\|\nabla\partial_t
b^\eps_R\|^2}{\eps_2}\|\partial_xb^\eps_R\|^2_{L^2}+C\|\nabla\partial_tb^\eps_R\|^2_{L^2}+C\|\nabla\partial_xb^\eps_R\|^2_{L^2}+C\frac{\|\nabla
u^\eps_R\|^2}{\eps_2}\|\partial_t\partial_xb^\eps_R\|^2_{L^2}\nonumber\\
&&+C\frac{\|\nabla\partial_x
u^\eps_R\|^2}{\eps_2}\|\partial_tb^\eps_R\|^2_{L^2}+C\frac{\|\nabla\partial_t
u^\eps_R\|^2}{\eps_2}\|\partial_xb^\eps_R\|^2_{L^2}\nonumber\\
&&+C\|\nabla\partial_tb^\eps_R\|^2_{L^2}+C\|\nabla\partial_xb^\eps_R\|^2_{L^2}+C\|\nabla
u^\eps_R\|^2_{L^2}+C\|\nabla
b^\eps_R\|^2_{L^2}.\label{est-pt10-2}\er Hence, we have \ber
&&\frac{d}{dt}\int
(|\partial_t\partial_xu^\eps_R|^2+|\partial_t\partial_xu^\eps_R|^2)+(1-\delta)\eps_1\int|\nabla\partial_t\partial_xu^\eps_R|^2
+(1-\delta)\eps_2\int|\nabla\partial_t\partial_xb^\eps_R|^2\nonumber\\
&\le&C\int(|\partial_t\partial_xu^\eps_R|^2+|\partial_t\partial_xb^\eps_R|^2)
+C(\frac{(\eps_1-\eps_2)^2}{\eps(\eps_1+\eps_2)}+\delta\eps)\int_0^t\int(|\nabla\partial_t\partial_xu^\eps_R|^2+|\nabla\partial_t\partial_xb^\eps_R|^2)\nonumber\\
&&+C\frac{\|\nabla
u^\eps_R\|^2}{\eps_1}\|\partial_t\partial_xu^\eps_R\|^2_{L^2}+C\frac{\|\nabla\partial_x
u^\eps_R\|^2}{\eps_1}\|\partial_tu^\eps_R\|^2_{L^2}+C\frac{\|\nabla\partial_t
u^\eps_R\|^2}{\eps_1}\|\partial_xu^\eps_R\|^2_{L^2}\nonumber\\
&&+C\frac{\|\nabla
b^\eps_R\|^2}{\eps_1}\|\partial_t\partial_xb^\eps_R\|^2_{L^2}+C\frac{\|\nabla\partial_x
b^\eps_R\|^2}{\eps_1}\|\partial_tb^\eps_R\|^2_{L^2}+C\frac{\|\nabla\partial_t
b^\eps_R\|^2}{\eps_1}\|\partial_xb^\eps_R\|^2_{L^2}\nonumber\\
&&+C\frac{\|\nabla
b^\eps_R\|^2}{\eps_2}\|\partial_t\partial_xu^\eps_R\|^2_{L^2}+C\frac{\|\nabla\partial_x
b^\eps_R\|^2}{\eps_2}\|\partial_tu^\eps_R\|^2_{L^2}+C\frac{\|\nabla\partial_t
b^\eps_R\|^2}{\eps_2}\|\partial_xb^\eps_R\|^2_{L^2}\nonumber\\
&&+C\frac{\|\nabla
u^\eps_R\|^2}{\eps_2}\|\partial_t\partial_xb^\eps_R\|^2_{L^2}+C\frac{\|\nabla\partial_x
u^\eps_R\|^2}{\eps_2}\|\partial_tb^\eps_R\|^2_{L^2}+C\frac{\|\nabla\partial_t
u^\eps_R\|^2}{\eps_2}\|\partial_xb^\eps_R\|^2_{L^2}\nonumber\\
&&+C\|\nabla\partial_tu^\eps_R\|^2_{L^2}+C\|\nabla\partial_xu^\eps_R\|^2_{L^2}
+C\|\nabla\partial_tb^\eps_R\|^2_{L^2}+C\|\nabla\partial_xb^\eps_R\|^2_{L^2}\nonumber\\
&&+C\|\nabla u^\eps_R\|^2_{L^2}+C\|\nabla b^\eps_R\|^2_{L^2}\nonumber\\
&&+C\frac{\beta_3(\eps)}\eps
+C\eps^{\kappa-1}+C\frac{(\eps_1-\eps_2)^2}{\eps(\eps_1+\eps_2)\min\{\eps_1,\eps_2\}}\beta_2(\eps)+C\frac{(\eps_1-\eps_2)^2}{\eps\sqrt\eps(\eps_1+\eps_2)}
\nonumber\\
&\le&C\int(|\partial_t\partial_xu^\eps_R|^2+|\partial_t\partial_xb^\eps_R|^2)
+C\eps_1\int_0^t\int|\nabla\partial_t\partial_xu^\eps_R|^2+C\eps_2\int_0^t\int|\nabla\partial_t\partial_xb^\eps_R|^2\nonumber\\
&&+C\frac{\|\nabla
u^\eps_R\|^2}{\eps_1}\|\partial_t\partial_xu^\eps_R\|^2_{L^2}+C\frac{\|\nabla
b^\eps_R\|^2}{\eps_1}\|\partial_t\partial_xb^\eps_R\|^2_{L^2}+C\frac{\|\nabla
b^\eps_R\|^2}{\eps_2}\|\partial_t\partial_xu^\eps_R\|^2_{L^2}+C\frac{\|\nabla
u^\eps_R\|^2}{\eps_2}\|\partial_t\partial_xb^\eps_R\|^2_{L^2}\nonumber\\
&&+C(\|\nabla\partial_x u^\eps_R\|^2+\|\nabla\partial_x
b^\eps_R\|^2)(\frac1{\eps_1}+\frac1{\eps_2})\beta_1(\eps)+C(\|\nabla\partial_t
u^\eps_R\|^2+\|\nabla\partial_t
b^\eps_R\|^2)(\frac1{\eps_1}+\frac1{\eps_2})\beta_2(\eps)\nonumber\\
&&+C\|\nabla\partial_tu^\eps_R\|^2_{L^2}+C\|\nabla\partial_xu^\eps_R\|^2_{L^2}
+C\|\nabla\partial_tb^\eps_R\|^2_{L^2}+C\|\nabla\partial_xb^\eps_R\|^2_{L^2}\nonumber\\
&&+C\|\nabla u^\eps_R\|^2_{L^2}+C\|\nabla b^\eps_R\|^2_{L^2}\nonumber\\
&&+C\frac{\beta_3(\eps)}\eps
+C\eps^{\kappa-1}+C\frac{(\eps_1-\eps_2)^2}{\eps(\eps_1+\eps_2)\min\{\eps_1,\eps_2\}}\beta_2(\eps)+C\frac{(\eps_1-\eps_2)^2}{\eps\sqrt\eps(\eps_1+\eps_2)}.
\label{est-pt10-3}\er Using Gronwall's inequality into
\eqref{est-pt10-3}, one gets that \ber
\int(|\partial_t\partial_xu^\eps_R|^2+|\partial_t\partial_xb^\eps_R|^2)
+\eps_1\int_0^t\int|\nabla\partial_t\partial_xu^\eps_R|^2+\eps_2\int_0^t\int|\nabla\partial_t\partial_xb^\eps_R|^2\le
C\beta_4(\eps).\label{est-pt-11}\er Here
\ber\beta_4(\eps)&=&\eps^\kappa+\frac{\beta_3(\eps)}\eps
+\eps^{\kappa-1}+C\frac{(\eps_1-\eps_2)^2}{\eps(\eps_1+\eps_2)\min\{\eps_1,\eps_2\}}\beta_2(\eps)+\frac{(\eps_1-\eps_2)^2}{\eps\sqrt\eps(\eps_1+\eps_2)}\nonumber\\
&&+\frac1{\min\{\eps_1,\eps_2\}}(\frac1{\eps_1}+\frac1{\eps_2})\beta_1(\eps)\beta_2(\eps)
+(\beta_1(\eps)+\beta_2(\eps))(\frac1{\eps_1}+\frac1{\eps_2})\nonumber\\
&&+\frac{\beta_0(\eps)}{\min\{\eps_1,\eps_2\}}.\label{beta4}\er Here
we require
$\frac{\beta_0(\eps)+\beta_1(\eps)+\beta_2(\eps)}{\min\{\eps_1,\eps_2\}}(\frac
1{\eps_1}+\frac 1{\eps_2})\le C.$

Hence
\ber&&\int|\partial_t\partial_x(u^\eps_R-b_R^\eps)|^2+(\eps_1+\eps_2)\int_0^t\int|\nabla\partial_t\partial_x(u^\eps_R-b^\eps_R)|^2\nonumber\\
&\le&
C(\frac{(\eps_1-\eps_2)^2}{(\eps_1+\eps_2)\min\{\eps_1,\eps_2\}}+\frac{\eps^2}{\min\{\eps_1,\eps_2\}})\beta_4(\eps)+\beta_3(\eps).\label{est-pt12}
\er Now we apply the anisotropic Sobolev imbedding inequality
\cite{TW2} and we get \ber \|(u^\eps_R,
b^\eps_R)\|_{L^\infty(\Omega\times (0, T))}&\le& C(\|(u^\eps_R,
b^\eps_R)\|_{L^\infty(0, T; H)}^{\frac12}\|\partial_x(u^\eps_R,
b^\eps_R)\|_{L^\infty(0, T;
L^2)}^{\frac 12}\nonumber\\
&&+\|(u^\eps_R, b^\eps_R)\|_{L^\infty(0, T;
L^2)}^{\frac12}\|\partial_x\partial_z(u^\eps_R,
b^\eps_R)\|_{L^\infty(0, T; L^2)}^{\frac 12})\nonumber\\
&\le&
C(\frac{\beta_1(\eps)}{\min\{\eps_1,\eps_2\}})^{\frac14}(\beta_2(\eps))^{\frac14}+C(\beta_0(\eps))^{\frac
14}(\frac{\beta_4(\eps)}{\min\{\eps_1,\eps_2\}})^{\frac 14}\to 0
\mbox{ when } \eps\to 0.\nonumber \er This completes the proof of
estimates \eqref{equ} and \eqref{res1c} in Theorem \ref{th1}.

\subsection{The Proof of Theorem \ref{thm1.2}}
Let $(\ue, \be)$ be the Leray-Hopf weak solutions to MHD systems
\eqref{vmhde1}-\eqref{vmhde5}. Decompose the solution as $(\ue,
\be)=(u^{\eps_1,0}+\urr, b^{\eps_1,0}+\bbee+\brr)$. Taking
$\nu_2^*=(\theta\eps_2)^{1+\tau}$ with $\tau\in [0, 1)$ and using
the system \eqref{hvmhd1}-\eqref{hvmhd5}, we have \ber
&&\pt\urr+\ue\cdot\nabla\urr+\urr\cdot\nabla u^{\eps_1,0}
-\eps_1\Delta\urr\nonumber\\
&=&-\nabla
(\pp-p^{\eps_1,0})+\be\cdot\nabla\brr+  b^{\eps_1,0}\cdot\nabla\bbee+\bbee\cdot\nabla\bbee+\brr\cdot\nabla\bbe\nonumber\\
&&+\bbee\cdot\nabla  b^{\eps_1,0}+\brr\cdot\nabla  b^{\eps_1,0}, \quad \mbox{ in }\quad \Omega\times(0, T),\label{mhdr1}\\
&&\pt\bbee+\pt\brr+\ue\cdot\nabla\brr+u^{\eps_1,0}\cdot\nabla\bbee+\urr\cdot\nabla\bbe
\nonumber\\
&&+\urr\cdot\nabla  b^{\eps_1,0} -\eps_2\Delta\brr-\eps_2\Delta
b^{\eps_1,0}-\eps_2\Delta\bbee\nonumber\\
&=&\be\cdot\nabla\urr+\bbee\cdot\nabla u^{\eps_1,0}+\brr\cdot\nabla u^{\eps_1,0}, \quad \mbox{ in }\quad \Omega\times(0, T),\label{mhdr2}\\
&&div\ue=div\be=div\urr=div\brr=0,\mbox{ in } \Omega\times(0, T),\label{mhdr3}\\
&&\urr=\brr=0,  \quad\mbox { on } \quad (x, y)\in \omega, z=0,h,
\quad 0\le t\le
T\label{mhdr4}\\
&&\urr(t=0)=\ue(0)-u^{\eps_1,0}(0),\nonumber\\
&&\brr(t=0)=\be(0)-  b^{\eps_1,0}(0)-\bbee(t=0), \quad \mbox{ on }
\quad \Omega .\label{mhdr5}\er By taking the scalar product of
\eqref{mhdr1} with $\urr$ and the scalar product of \eqref{mhdr2}
with $\brr$, we have \ber &&\frac 12\frac
d{dt}\int(|\urr|^2+|\brr|^2)+\eps_1\int|\nabla\urr|^2+\eps_2\int|\nabla\brr|^2
=\sum_{i=1}^{11}I_i,\label{est-21} \er where $I_i, i=1, \cdots, 11,$
are given respectively as follows \ber
&&I_1=-\int\nabla(\pp-\ppp)\urr;\quad I_2=\int\pt\bbee\brr;\nonumber\\
&&I_3=-\int\ue\cdot\nabla\urr\urr-\int\ue\cdot\nabla\brr\brr;\nonumber\\
&&I_4=\int  b^{\eps_1,0}\cdot\nabla\bbee\urr-\int
u^{\eps_1,0}\cdot\nabla\bbe\brr
\nonumber\\
&&I_5=\int\bbee\cdot\nabla\bbee\urr;\nonumber\\
&&I_6=-\int\urr\cdot\nabla\bbee\brr+\int\brr\cdot\nabla\bbee\urr\nonumber\\
&&I_7=\int\bbee\cdot\nabla  b^{\eps_1,0}\urr+\int\bbee\cdot\nabla u^{\eps_1,0}\brr;\nonumber\\
&&I_8=-\int\urr\cdot\nabla u^{\eps_1,0}\urr-\int\urr\cdot\nabla  b^{\eps_1,0}\brr\nonumber\\
&&\quad +\int\brr\cdot\nabla  b^{\eps_1,0}\urr+\int\brr\cdot\nabla u^{\eps_1,0}\brr; \nonumber\\
&&I_{9}=\int\be\cdot\nabla\brr\urr+\int\be\cdot\nabla\urr\brr;\nonumber\\
&&I_{10}=\eps_2\int\Delta b^{\eps_1,0}\brr;\quad
I_{11}=\eps_2\int\Delta  \bbee\brr; \nonumber\er We now bound each
of $I_i, i=1, \cdots, 11$.

First, similar to the estimates of $J_1, \cdots, J_5, J_7,\cdots
J_{13}$, we can estimate $I_1, \cdots, I_5, I_7,\cdots, I_{11}$ to
get \ber I_1+\cdots+I_5+I_7+\cdots+I_{11}\le
C\int(|u_R^\eps|^2+|b_R^\eps|^2)+C(\frac{(\sqrt{\eps_2})^{1-\tau}}{\sqrt\theta}+(\sqrt{\theta\eps_2})^{1+\tau}).\label{est-22}\er
Secondly, we estimate $I_6$ as follows: \ber
I_6=I_{61}+I_{62},\nonumber\er where \ber
I_{61}=-\int\urr\cdot\nabla\bbe\brr ;\quad
I_{62}&=&\int\brr\cdot\nabla\bbe\urr .\nonumber \er Now we estimate
$I_{61}$, which is split into four parts: \ber I_{61}&=&-\int
U_R^{\eps}\cdot\nabla_{x,y}B_B^{\eps}B_R^{\eps}-\int
U_R^{\eps}\cdot\nabla_{x,y}\bbet\brrt\nonumber\\
&&-\int\urrt\partial_z\bbet\brrt
-\int\urrt\partial_zB_B^{\eps}B_R^{\eps}\nonumber\\
&=&K_1+K_2+K_3+K_4.\nonumber\er The integrals $K_1, K_2, K_3$ can be
easily bounded by \ber K_1+K_2+K_3\le C\|(u^{\eps_1,0},
  b^{\eps_1,0})\|_{H^s}\|(\urr, \brr)\|_{L^2}^2, \quad s>\frac52.\nonumber\er For $K_4$, we have \ber
K_4&=&-\int_\omega\int_0^{\frac
h4}\urrt\cdot\partial_zB_{B+}^{\eps}B_R^{\eps}-\int_\omega\int_{\frac
{3h}4}^h\urrt\cdot\partial_zB_{B-}^{\eps}B_R^{\eps}\nonumber\\
&=&-\int_\omega\int_0^{\frac h4}\frac{\urrt}z
z^2\partial_zB_{B+}^{\eps}\frac{B_R^{\eps}}z-\int_\omega\int_{\frac
{3h}4}^h\frac{\urrt}{h-z}(h-z)^2\partial_zB_{B-}^{\eps}\frac{B_R^{\eps}}{h-z}\nonumber\\
&\le&\|\frac{\urrt}z\|_{L^2}\|z^2\partial_zB_{B+}^{\eps}\|_{L^\infty}\|\frac{B_R^{\eps}}z\|_{L^2}
+\|\frac{\urrt}{h-z}\|_{L^2}\|(h-z)^2\partial_zB_{B-}^{\eps}\|_{L^\infty}
\|\frac{B_R^{\eps}}{h-z}\|_{L^2}\nonumber\\
&\le&C\|(u^{\eps_1,0}, b^{\eps_1,0})\|_{H^s}\sqrt{(\theta\eps_2)^{1+\tau}}\|\partial_z{\urrt}\|_{L^2}\|\partial_zB_R^{\eps}\|_{L^2}\nonumber\\
&\le&\theta^{1+\tau}(\eps_2)^{1+\frac\tau2}\|\partial_z\brr\|_{L^2}^2+C\theta^{1+\tau}(\eps_2)^{\frac\tau2}\|(u^{\eps_1,0},
b^{\eps_1,0})\|_{H^s}^2
\|\partial_zu_R^{\eps}\|_{L^2}^2\nonumber\\
&\le&\theta^{1+\tau}\eps_2\|\partial_z\brr\|_{L^2}^2+C\theta^{1+\tau}\|\partial_zu_R^{\eps}\|_{L^2}^2\nonumber\er
for some $\theta>0$ sufficiently small if $\eps_2\to 0$. Here we
used the Hardy's inequality $\|(\frac{f(z)}z,
\frac{f(z)}{h-z})\|_{L^2(0, 1)}\le C\|\partial_z f(z)\|_{L^2(0, 1)}$
when $f(0)=0$.

Hence, we have \ber I_{61}\le
\theta^{1+\tau}\eps_2\|\partial_z\brr\|_{L^2}^2+C\theta^{1+\tau}\|\partial_zu_R^{\eps}\|_{L^2}^2+C(\|\urr|_{L^2}^2+
\|\brr\|_{L^2}^2).\nonumber\er Similarly, $I_{62}$ can be bounded by
\ber I_{62}\le
\theta^{1+\tau}\eps_2\|\partial_z\brr\|_{L^2}^2+C\theta^{1+\tau}\|\partial_zu_R^{\eps}\|_{L^2}^2+C(\|\urr|_{L^2}^2+
\|\brr\|_{L^2}^2).\nonumber\er Thus, we have the estimate on
$I_6$\ber
I_6\le\theta^{1+\tau}\eps_2\|\partial_z\brr\|_{L^2}^2+C\theta^{1+\tau}\|\partial_zu_R^{\eps}\|_{L^2}^2+C(\|\urr|_{L^2}^2+
\|\brr\|_{L^2}^2).\label{est-23}.\er Putting
(\ref{est-22})-(\ref{est-23}) into (\ref{est-21}), one have \ber
&&\frac 12\frac
d{dt}\int(|\urr|^2+|\brr|^2)+(\eps_1-C\theta^{1+\tau})\int|\nabla\urr|^2+\eps_2(1-\theta^{1+\tau})\int|\nabla\brr|^2\nonumber\\
&\le&C\int(|u_R^\eps|^2+|b_R^\eps|^2)+C(\frac{(\sqrt{\eps_2})^{1-\tau}}{\sqrt\theta}+(\sqrt{\theta\eps_2})^{1+\tau}).\label{est-24}\er
Taking $\theta>0$ to be sufficiently small and to be independent of
$\eps_2$ and applying Gronwall's inequality to \eqref{est-24} and
using the assumption \eqref{initialass1} on initial data, we can get
the estimate rate \eqref{estrate1} in Theorem \ref{thm1.2}.

\vspace{0.3cm} \noindent{\bf Acknowledgments} Wang's Research is
supported by NSFC(No.11371042). Xin's research is also supported
partially Zhang Ge Ru Foundation, Hong Kong RDG Earmarked Research
Grants CUHK-14305315 and CUHK 4048/13P, NSFC/RDG Joint Research
Scheme Grant N-CUHK 443-14, and a Focus Area Grant from the Chinese
University of Hong Kong.

\end{document}